\documentclass[12pt,a4paper,reqno]{amsart}
\usepackage[english]{babel}
\usepackage{amssymb,latexsym,amsfonts,amsthm,upref,amsmath}
\usepackage[margin=1in]{geometry}
\usepackage[dvipsnames,x11names]{xcolor}
\definecolor{myblue}{HTML}{003399}
\usepackage{hyperref}
\hypersetup{colorlinks,citecolor=myblue,filecolor=black,linkcolor=myblue,urlcolor=myblue}
\usepackage{enumerate} 
\usepackage{tikz}
\usepackage{float}
\usepackage{cleveref}
\usepackage{mathtools}
\usepackage{cite}
\usepackage{filecontents}
\makeatletter
\newcommand{\leqnomode}{\tagsleft@true}
\newcommand{\reqnomode}{\tagsleft@false}
\newcommand{\cev}[1]{\reflectbox{\ensuremath{\vec{\reflectbox{\ensuremath{#1}}}}}}
\makeatother

\newtheorem*{corintro*}{Corollary}
\newtheorem*{thm*}{Theorem}
\newtheorem*{lem*}{Lemma}
\newtheoremstyle{prim}{}{}{\normalfont}{}{\bfseries}{.}{ }{}
\newtheoremstyle{stil}{}{}{\slshape}{}{\bfseries}{.}{ }{}
\theoremstyle{stil}
\newtheorem{thm}{Theorem}[section]
\newtheoremstyle{defi}{}{}{}{}{\bfseries}{.}{ }{}
\theoremstyle{defi}
\newtheorem{defn}[thm]{Definition}
\theoremstyle{defi}
\newtheorem{rem}[thm]{Remark}
\theoremstyle{stil}
\newtheorem{pro}[thm]{Proposition}
\theoremstyle{stil}
\newtheorem{lem}[thm]{Lemma}
\theoremstyle{stil}
\newtheorem{kor}[thm]{Corollary}
\theoremstyle{prim}
\newtheorem{ex}[thm]{Example}
\newenvironment{prf}{\noindent \textit{Proof.}}{\null\hfill$\qed$\hskip
2mm\vskip 2mm}

\newcommand{\sgn}{ {\rm sgn}\ts}

\newcommand{\vac}{\mathop{\mathrm{\boldsymbol{1}}}}
\newcommand{\ndo}{\mathop{\mathrm{End}}}
\newcommand{\om}{\mathop{\mathrm{Hom}}}

\newcommand{\Sym}{\mathfrak S}

\newcommand{\diag}{\mathop{\mathrm{diag}}}

\newcommand{\DY}{ {\rm DY}}
\newcommand{\Y}{ {\rm Y}}
\newcommand{\B}{ {\rm B}}
\newcommand{\gl}{\mathfrak{gl}}
\newcommand{\sll}{\mathfrak{sl}}

\newcommand{\hsym}{\mathop{\hspace{2pt}\underset{h^n}{\sim}\hspace{2pt}}}

\newcommand{\wtld}{\widetilde}
\newcommand{\wht}{\widehat}
\newcommand{\wvr}{\overline}
\newcommand{\wndr}{\underline}

\newcommand{\ot}{\otimes}

\newcommand{\ts}{\,}
\newcommand{\tss}{\hspace{1pt}}

\newcommand{\CC}{\mathbb{C}\tss}

\newcommand{\ZZ}{\mathbb{Z}\tss}
\newcommand{\TT}{\mathbb{T}}
\newcommand{\BB}{\mathbb{A}}

\newcommand{\XX}{\mathbb{B}}

\newcommand{\Ac}{ {\rm A}}

\newcommand{\Sc}{\mathcal{S}}
\newcommand{\Tc}{\mathcal{T}}
\newcommand{\Bc}{\mathcal{B}}
\newcommand{\Ec}{\mathcal{E}}

\newcommand{\Uc}{\mathcal{U}}
\newcommand{\Vc}{\mathcal{V}}
\newcommand{\Wc}{\mathcal{W}}

\newcommand{\g}{\mathfrak{g}}

\newcommand{\z}{\mathfrak{z}}

\newcommand{\sdet}{ {\rm sdet}\ts}

\newcommand{\qdet}{ {\rm qdet}\ts}
\newcommand{\tr}{ {\rm tr}}
\newcommand{\fand}{\quad\text{and}\quad}
\newcommand{\Fand}{\qquad\text{and}\qquad}

\newcommand{\non}{\nonumber}
\newcommand{\beq}{\begin{equation}}
\newcommand{\eeq}{\end{equation}}
\newcommand{\ben}{\begin{equation*}}
\newcommand{\een}{\end{equation*}}

\begin{document}

\title[Quasi modules for the quantum affine vertex algebra in type $A$]{Quasi  modules for the quantum affine vertex algebra in type $A$}

\author{Slaven Ko\v{z}i\'{c}} 
\address{School of Mathematics and Statistics F07, University of Sydney, NSW 2006, Australia}
\address{Department of Mathematics, Faculty of Science, University of Zagreb, 10000 Zagreb, Croatia}
\email{kslaven@math.hr}

\maketitle

\begin{abstract}
We consider the quantum affine vertex algebra $\mathcal{V}_{c}(\mathfrak{gl}_N)$ associated with the rational $R$-matrix, as defined by Etingof and Kazhdan. We introduce certain subalgebras $\textrm{A}_c (\mathfrak{gl}_N)$ of the completed double Yangian $\widetilde{\textrm{DY}}_{c}(\mathfrak{gl}_N)$ at the level $c\in\mathbb{C}$, associated with the reflection equation, and we employ their structure to construct examples of quasi $\mathcal{V}_{c}(\mathfrak{gl}_N)$-modules. Finally, we use the quasi module map, together with the explicit description of the center of $\mathcal{V}_{c}(\mathfrak{gl}_N)$, to obtain  formulae for families of central elements in the completed algebra $\widetilde{\textrm{A}}_c (\mathfrak{gl}_N)$.
\end{abstract}

\section*{Introduction}

 In order to describe  integrable systems with the boundary conditions,  E. K. Sklyanin introduced in \cite{S} the {\em reflection algebras}, a   class of algebras associated with   $R$-matrix $R(u)$ which are defined by the {\em reflection equation}
\beq\label{reflectionequation}
R_{12}(u-v)B_1 (u)R_{12}(u+v)B_2 (v)=B_2 (v)R_{12}(u+v)B_1 (u)R_{12}(u-v).
\eeq
We explain the precise meaning of \eqref{reflectionequation}  in Section \ref{explinedin}.
His approach was motivated by   Cherednik's 
 treatment of factorized scattering with reflection  \cite{C}.
Furthermore, Sklyanin constructed an analogue of the {\em quantum determinant}    and described the {\em algebraic Bethe ansatz}; see \cite{S}. Later on, different classes of algebras  defined via relations of the form similar to or same as \eqref{reflectionequation} were extensively studied; see, e.g., \cite{KS,KJC,MRS,MR,RS}.

In this paper, we consider certain family of reflection algebras associated with Yang  $R$-matrix,  studied by A. Molev and E. Ragoucy in \cite{MR}, which are coideal subalgebras of the Yangian $\Y(\gl_N)$. We introduce the subalgebra $\Ac_c (\gl_N)$ of the $h$-adically completed double Yangian $\wtld{\DY}_c (\gl_N)$ at the level $c\in\CC$ which, roughly speaking,  consists  of   two reflection algebras.
Motivated by the correspondence,   indicated in \cite{EK}, between the {\em $\mathcal{S}$-locality} (see \eqref{locality} below) and the reflection equation which appeared in  work of N. Yu. Reshetikhin and M. A. Semenov-Tian-Shansky  \cite{RS}, we investigate  algebras $\Ac_c (\gl_N)$ using the theory of quantum VOAs.

The notion of {\em quantum vertex operator algebra}
 ({\em quantum VOA}) was introduced by P. Etingof and D. Kazhdan in \cite{EK}. {\em Quantum affine VOA} can be  associated with rational, trigonometric and elliptic $R$-matrix; see \cite{EK}. In the rational case,  the   double Yangian $\DY(\gl_N)$ over $\CC[[h]]$ can be used to define the quantum VOA structure   on its vacuum module $\Vc_c(\gl_N)$ at the level $c\in\CC$. The theory of quantum vertex algebras was further developed and generalized by H.-S. Li; see, e.g., \cite{Li4,Li} and references therein. In particular, certain more general objects, such as {\em $h$-adic nonlocal vertex algebras} and their {\em  quasi modules}, were introduced and  studied in \cite{Li}. The main result of this paper  is a construction of   the quasi module map 
$Y_{\Wc_c(\gl_N)}$ on $\Vc_{2c}(\gl_N)$,
 so that  the {\em vacuum module} $\Wc_c(\gl_N)$ for the algebra $\Ac_c (\gl_N)$ acquires a  quasi $\Vc_{2c}(\gl_N)$-module structure.

We use the quasi module map to obtain further information on the algebra $\Ac_c (\gl_N)$.
In our previous paper \cite{JKMY}, coauthored with N. Jing,   Molev and F. Yang,
the center $\z(\Vc_c(\gl_N))$ of the quantum VOA  $\Vc_c(\gl_N)$  was described  by providing explicit formulae for its algebraically independent topological generators, thus establishing the quantum analogue of the  Feigin--Frenkel theorem in type $A$; see \cite{FF,CT,CM}. By considering the image of the center $\z(\Vc_{-N}(\gl_N))$, with respect to  the quasi module map $Y_{\Wc_{-N/2}(\gl_N)}$, we find explicit formulae for families of central elements   in the completed  algebra $\wtld{\Ac}_{-N/2} (\gl_N)$, which are, due to the {\em fusion procedure}  originated in the
work of  A. Jucys \cite{J},  parametrized  by arbitrary partitions with at most $N$ parts. For $c\neq -N$ we obtain only one family of central elements in $\wtld{\Ac}_{c/2} (\gl_N)$, which, roughly speaking, coincide  with the  coefficients of the product of two {\em Sklyanin determinants} (i.e. with the  coefficients of the product of four quantum determinants); see \cite{S,MR}.
In the end, we employ these central elements to obtain invariants of the vacuum module $\Wc_c(\gl_N)$.

\section{Reflection algebras}\label{sec1}
\numberwithin{equation}{section}
In this section, we  recall the definition of the double Yangian $\DY(\gl_N)$ over $\CC[[h]]$; see \cite{ I}. Next, we follow \cite{MR} to introduce a certain class of reflection algebras.  We employ their structure   to define subalgebra $\Ac_c(\gl_N)$ of the $h$-adically completed double Yangian $\wtld{\DY}_c (\gl_N)$ at the level $c\in\CC$,   which will be our main point of interest in this paper. 

\subsection{Double Yangian for \texorpdfstring{$\gl_N$}{glN}}\label{secc1}
Let $N\geqslant 2$ be an integer and let $h$ be a formal parameter.
Denote by $R(u)$  the {\em Yang $R$-matrix} over $\CC[[h]]$ defined by
\beq\label{yang}
R(u)=1-hPu^{-1},
\eeq
where $1$ is the identity and $P$ is the permutation operator in $\CC^N\ot \CC^N $, $P\colon x\ot y\mapsto y\ot x$. 
$R$-matrix \eqref{yang} satisfies the {\em Yang--Baxter equation}
\beq\label{ybe}
R_{12}(u)\ts R_{13}(u+v)\ts R_{23}(v)
=R_{23}(v)\ts R_{13}(u+v)\ts R_{12}(u).
\eeq
Both sides of \eqref{ybe} are operators on the triple tensor product $(\CC^{N})^{\ot 3}$ and subscripts indicate the copies of $\CC^N$ on which $R(u)$ acts, for example, 
$R_{12}(u)=R(u)\ot 1$ and $R_{23}(v)=1\ot R (v)$.
Let $g(u)$  be
the unique   series in $1+u^{-1}\CC[[u^{-1}]]$
satisfying
\beq\label{geq}g(u+N)=g(u)(1-u^{-2}).\eeq 
The $R$-matrix $\overline{R}(u)=\wvr{R}_{12}(u)=g(u/h)R(u)$ possesses the {\em crossing symmetry} properties,
\beq\label{csym}
\left(\wvr{R}_{12}(u)^{-1}\right)^{t_1}\wvr{R}_{12}(u+hN)^{t_1}=1\fand \left(\wvr{R}_{12}(u)^{-1}\right)^{t_2}\wvr{R}_{12}(u+hN)^{t_2}=1,
\eeq
where
$t_i$ denotes the transposition applied on the tensor factor $i=1,2$;
and the {\em unitarity property}
\beq\label{uni}
\wvr{R}_{12}(u)\wvr{R}_{12}(-u)=1,
\eeq
see, e.g., \cite[Section 2]{JKMY} for more details.

\allowdisplaybreaks

The {\em double Yangian} $\DY(\gl_N)$ for $\gl_N$  is defined as the associative algebra over $\CC[[h]]$ generated by the   central element $C$ and the elements
$t_{ij}^{(\pm r)}$, where $i,j=1,\ldots N$ and $r=1,2,\ldots$,    subject to the following defining relations (see \cite{I}),  
\begin{align}
R\big(u-v \big)\ts T_1(u)\ts T_2(v)&=T_2(v)\ts T_1(u)\ts
R\big(u-v \big),\label{RTT2}\\
R\big(u-v \big)\ts T^+_1(u)\ts T^+_2(v)&=T^+_2(v)\ts T^+_1(u)\ts
R\big(u-v \big),\label{RTT1}\\
\wvr{R}\big(u-v+hC/2\big)\ts T_1(u)\ts T^+_2(v)&=T^+_2(v)\ts T_1(u)\ts
\wvr{R}\big(u-v-hC/2\big).\label{RTT3}
\end{align}
The elements $T(u)$ and $T^+(u)$ in $\ndo\CC^N \ot \DY(\gl_N)[[u^{\pm 1}]]$ are defined by
\begin{align*}
T(u)=\sum_{i,j=1}^N e_{ij}\ot t_{ij}(u)\Fand T^{+}(u)=\sum_{i,j=1}^N e_{ij}\ot t_{ij}^{+}(u),
\end{align*}
where the $e_{ij}$ are the matrix units,
and the series $t_{ij}(u) $ and $t_{ij}^+ (u) $ are given by
$$t_{ij}(u)=\delta_{ij}+h\sum_{r=1}^{\infty} t_{ij}^{(r)}u^{-r}\Fand t_{ij}^+ (u)=\delta_{ij}-h\sum_{r=1}^{\infty} t_{ij}^{(-r)}u^{r-1}.$$
We use the  subscript to indicate a copy of the matrix in the tensor product algebra
$ (\ndo\CC^N)^{\ot m}\ot\DY(\gl_N)$,
so that, for example,
\beq\label{subscript}
T_k(u)=\sum_{i,j=1}^{N} 1^{\otimes (k-1)} \ot e_{ij} \ot 1^{\ot (m-k)} \ot t_{ij}(u).
\eeq
In particular, we have $m=2$ in defining relations \eqref{RTT2}--\eqref{RTT3}.

The {\em Yangian} $\Y(\gl_N)$ is the subalgebra of  $\DY(\gl_N)$ generated by the elements $t_{ij}^{(r)}$, $i,j=1,\ldots,N$, $r=1,2,\ldots$.  The {\em dual Yangian } $\Y^+ (\gl_N)$ is the subalgebra of the double Yangian $\DY(\gl_N)$ generated by the elements $t_{ij}^{(-r)}$, $i,j=1,\ldots,N$, $r=1,2,\ldots$. 
For any complex number $c$ denote by $\DY_c (\gl_N)$ the {\em double Yangian at the level $c$}, i.e. the quotient of the algebra $\DY(\gl_N)$ by the ideal generated by the element $C-c$.

Recall that the {\em $h$-adic topology} on an arbitrary $\CC[[h]]$-module $V$ is the topology generated by the  basis $v+h^n V$, $v\in V$, $n\in\mathbb{Z}_{\geqslant 1}$. 
The {\em vacuum module $\Vc_c(\gl_N)$ at the level $c$} over the double Yangian is the $h$-adic completion of the quotient of the algebra $\DY_c (\gl_N)$ by the left ideal generated by all elements $t_{ij}^{(r)}$, $ r=1,2,\ldots$, i.e. the $h$-adic completion of
\beq\label{imageofone}
\DY_c(\gl_N) /\DY_c(\gl_N)  \big\langle t_{ij}^{(r)}\,:\,i,j=1,\ldots,N,\, r=1,2,\ldots\big\rangle. 
\eeq
By the Poincar\'e--Birkhoff--Witt theorem for the double Yangian, see \cite[Theorem 2.2]{JKMY},  the  vacuum module $\Vc_c(\gl_N)$ is isomorphic, as a $\CC[[h]]$-module, to the $h$-adically completed  dual Yangian $\widehat{\Y}^+(\gl_N)$.

\subsection{Algebra \texorpdfstring{$\Ac_{c}(\gl_N)$}{Ac(glN)}
}\label{explinedin}
We now proceed as in \cite{MR} to introduce the reflection algebras.
Fix nonnegative integer $M\leqslant N$. Let $G=(g_{ij})_{i,j=1}^N$ be the diagonal matrix of order $N$,
\beq\label{diagonal}
G=\diag(\varepsilon_1,\ldots,\varepsilon_N),
\eeq
where $\varepsilon_1=\ldots=\varepsilon_{M}=1$ and $\varepsilon_{M+1}=\ldots=\varepsilon_{N}=-1$.
Let $c$ be a fixed complex number. Consider the series 
\begin{align}
&B^{+}(u)=\sum_{i,j=1}^N e_{ij}\ot b_{ij}^{+}(u)\in \ndo\CC^N \ot\widehat{\Y}^+ (\gl_N)[[u]]\fand\label{B+}\\
&B(u)=\sum_{i,j=1}^N e_{ij}\ot b_{ij}(u)\in \ndo\CC^N \ot\Y (\gl_N)[[u^{-1}]]\label{B}
\end{align}
defined by
\beq\label{Bseries}
B^{+}(u)=T^{+}(u)\tss G\tss T^{+}(-u)^{-1}\Fand B(u)=T(u+hc)\tss G\tss T(-u)^{-1}.
\eeq
We can write the matrix entries of \eqref{B+} and \eqref{B} as
\begin{align*}
&b_{ij}^{+}(u)=g_{ij} -h\sum_{r=1}^\infty b_{ij}^{(-r)}u^{r-1}\Fand b_{ij}(u)=g_{ij} +h\sum_{r=1}^\infty b_{ij}^{(r)}u^{-r}
\end{align*}
for some elements  $b_{ij}^{(-r)}\in\widehat{\Y}^+(\gl_N)$ and $b_{ij}^{(r)}\in\Y(\gl_N)$.

Series \eqref{Bseries}   satisfy the {\em unitary condition}
\beq\label{rel1}
B^{+}(u)B^{+}(-u)=1\Fand B(u)B(-u-hc)=1.
\eeq
Furthermore, using \eqref{RTT2}--\eqref{RTT3} and
$R(u)G_1 R(v) G_2 =G_2 R(v) G_1 R(u)$
 one can easily verify that the following {\em reflection relations} hold for the elements of  the $h$-adically completed double Yangian  $\widehat{\DY}_c(\gl_N)$ at the level $c\in \CC$:
\begin{align}
&R(u-v)B_1^+(u) R(u+v) B_2^+(v)=B_2^+(v)R(u+v)B_1^+(u)R(u-v)\label{RBRB1},\\
&R(u-v)B_1(u) R(u+v+hc) B_2(v)=B_2(v)R(u+v+hc)B_1(u)R(u-v)\label{RBRB2},\\
&\wvr{R}(u-v+3hc/2)B_1(u) \wvr{R}(u+v-hc/2) B_2^+(v)\non\\
&\qquad\qquad\qquad\qquad\qquad\qquad=B_2^+(v)\wvr{R}(u+v+3hc/2)B_1(u)\wvr{R}(u-v-hc/2).\label{RBRB3}
\end{align}
As in Section \ref{secc1}, the subscripts in \eqref{RBRB1}--\eqref{RBRB3} indicate a copy of the matrix in the tensor product algebra $(\ndo\CC^N)^{\ot 2} \ot \widehat{\DY}_c (\gl_N)$; recall \eqref{subscript}.

For  $i,j=1,\ldots,N$ and $r=1,2,\ldots$ 
let $\Ac_c '(\gl_N)$ be  the subalgebra of   $\widehat{\DY}_c(\gl_N)$  generated by the elements $b_{ij}^{(-r)}$ and $b_{ij}^{(r)}$, 
let $\B'^+ (\gl_N)$ be the subalgebra of the $h$-adically completed dual Yangian $\widehat{\Y}^{+}(\gl_N)$ generated by the elements $b_{ij}^{(-r)}$ and
let $\B'_c  (\gl_N)$ be the subalgebra of the Yangian  $\Y(\gl_N)$ generated by the elements  $b_{ij}^{(r)}$.
\begin{rem}\label{remarkk}
By setting $h=1$ in the algebra $\B'_0  (\gl_N)$ we obtain  the reflection algebra
$\mathcal{B}(N,N-M)$ over $\CC$, as defined in \cite{MR}.
\end{rem}

As in \cite{Li}, for an arbitrary $\CC[[h]]$-submodule  $V$ of  $\widehat{\DY}_c(\gl_N)$ we define
\beq\label{xcvb}
[V] =\left\{v\in \widehat{\DY}_c(\gl_N) \,:\, h^n v\in V\text{ for some } n\geqslant 0  \right\}.
\eeq
Finally, consider the following subalgebras of $\widehat{\DY}_c(\gl_N)$:
$$\Ac_c (\gl_N)=[\Ac_c '(\gl_N)],\quad \B_c  (\gl_N)=[\B'_c  (\gl_N)]\fand \B^+ (\gl_N)=[\B'^+ (\gl_N)].$$
Clearly, the following inclusions hold:
$$\Ac_c (\gl_N)\subset \widehat{\DY}_c (\gl_N),\quad \B_c  (\gl_N)\subset\Y  (\gl_N) \fand\B^+ (\gl_N)\subset\widehat{\Y}^+ (\gl_N) .$$
Moreover, due to \cite[Lemma 3.5]{Li}, the induced  topology on $\Ac_c (\gl_N)$,  $\B_c  (\gl_N)$ and $\B^+ (\gl_N)$ from  $\widehat{\DY}_c (\gl_N)$ coincides with the $h$-adic topology on these algebras.

We  now introduce some new notation in order to write the more general form of relations \eqref{RBRB1}--\eqref{RBRB3}. For positive integers $n,m$ and the families of variables $u=(u_1,\ldots,u_n)$ and $v=(v_{1},\ldots, v_{m})$ set
$$\wvr{R}_{ij}=\wvr{R}_{ij}(u_i-v_{j-n})\fand \wndr{\wvr{R}}_{i j}=\wvr{R}_{ij}(u_i+v_{j-n}), \qquad  i=1,\ldots,n,\,j=n+1,\ldots,n+m.$$
 Consider the functions  with values in the space
$(\ndo\CC^N)^{\ot n} \ot (\ndo\CC^N)^{\ot m}$
\begin{align}
& \wvr{R}_{nm}^{12}(u|v)=\prod_{i=1,\ldots,n}^{\longrightarrow} \prod_{j=n+1,\ldots,n+m}^{\longleftarrow}\wvr{R}_{ij}\Fand
\wndr{\wvr{R}}_{nm}^{12}(u|v)=\prod_{i=1,\ldots,n}^{\longrightarrow} \prod_{j=n+1,\ldots,n+m}^{\longrightarrow}\wndr{\wvr{R}}_{ij}\label{rmatrixces}
\end{align}
with the arrows indicating the order of the factors.
The functions   $R_{nm}^{12}(u|v)$ and $\wndr{R}_{nm}^{12}(u|v)$   corresponding to  $R$-matrix \eqref{yang} can be defined analogously.
Introduce the series 
\begin{align}
&\wndr{B}_{n}^{+}(u)=\prod_{i=1,\ldots,n}^{\longrightarrow} \left(B_{i}^+ (u_i) \wvr{R}_{i\,i+1}(u_i+u_{i+1})\ldots \wvr{R}_{in}(u_i+u_{n})  \right)  \fand\label{B+u}\\
&\wndr{B}_{n}(u)=\prod_{i=1,\ldots,n}^{\longrightarrow} \left(B_{i} (u_i) \wvr{R}_{i\,i+1}(u_i+u_{i+1}+hc)\ldots \wvr{R}_{in}(u_i+u_{n}+hc)   \right).\label{Bu}
\end{align}

For a family of variables $u=(u_1,\ldots, u_n)$ and $\alpha\in\CC$ we will often denote the families $(u_1+\alpha h,\ldots, u_n+\alpha h)$ and $(\alpha u_1,\ldots ,\alpha u_n)$ by $u+\alpha h$ and $\alpha u$ respectively.
 We also adopt  the superscript notation
for multiple tensor products of the form
$$(\ndo\mathbb{C}^{N})^{\otimes n} \otimes
(\ndo\mathbb{C}^{N})^{\otimes m}\otimes
(\ndo\mathbb{C}^{N})^{\otimes k}\otimes \Ac_c(\gl_N)
\otimes \Ac_c(\gl_N)\otimes \Ac_c(\gl_N).
$$
Expressions like $\wndr{B}_{n}^{+14}(u)$ or $\wndr{B}_{k}^{ 35}(w)$, where $w=(w_1,\ldots ,w_k)$, will be understood
as the respective operators $\wndr{B}_{n}^+ (u)$ or $\wndr{B}_{k}(w)$, whose non-identity components
belong to the corresponding tensor factors. In particular, the non-identity
components of $\wndr{B}_{k}^{ 35}(w)$ belong to the
factors
$
n+m+1,\,n+m+2,\,\dots,\,n+m+k$  and $n+m+k+2
$.
This notation is employed in the next proposition, which can be proved using  \eqref{RBRB1}--\eqref{RBRB3} and Yang--Baxter equation \eqref{ybe}.
\begin{pro}\label{rbrb_pro}
For any positive integers $n$ and $m$ the following equalities hold on
$(\ndo\CC^N)^{\ot n}\ot(\ndo\CC^N)^{\ot m}\ot \Ac_c(\gl_N)$:
\begin{align}
&R_{nm}^{12}(u|v)\wndr{B}_{n}^{+13}(u)\wndr{R}_{nm}^{12}(u|v)\wndr{B}_{m}^{+23}(v)
=\wndr{B}_{m}^{+23}(v)\wndr{R}_{nm}^{12}(u|v)\wndr{B}_{n}^{+13}(u)R_{nm}^{12}(u|v),\label{rbrb1}\\
&R_{nm}^{12}(u|v)\wndr{B}_{n}^{13}(u)\wndr{R}_{nm}^{12}(u+hc|v)\wndr{B}_{m}^{23}(v)
=\wndr{B}_{m}^{23}(v)\wndr{R}_{nm}^{12}(u+hc|v)\wndr{B}_{n}^{13}(u)R_{nm}^{12}(u|v),\label{rbrb2}\\
&\wvr{R}_{nm}^{12}(u+3hc/2|v)\wndr{B}_{n}^{13}(u)\wndr{\wvr{R}}_{nm}^{12}(u-hc/2|v)\wndr{B}_{m}^{+23}(v)\non\\
&\qquad\qquad\qquad\qquad\qquad =\wndr{B}_{m}^{+23}(v)\wndr{\wvr{R}}_{nm}^{12}(u+3hc/2|v)\wndr{B}_{n}^{13}(u)\wvr{R}_{nm}^{12}(u-hc/2|v).\label{rbrb3}
\end{align}
\end{pro}

Our next goal is to derive \eqref{csym5}, which will be useful in what follows.
First, note that by applying the transposition $t_1$ on the first and $t_2$ on the second equality in    \eqref{csym} we get
\beq\label{csym2}
{}^{rl\hspace{-3pt}}\left(\wvr{R}_{12}(u)^{-1}\right)\cdot\wvr{R}_{12}(u+hN)=1
\Fand
{}^{lr\hspace{-3pt}}\left(\wvr{R}_{12}(u)^{-1}\right)\cdot\wvr{R}_{12}(u+hN)=1,
\eeq
where the superscript $ rl $ ($ lr $) in \eqref{csym2} indicates that the first tensor factor of $\wvr{R}_{12}(u)^{-1}$ is applied from the right (left) while the second tensor factor  of $\wvr{R}_{12}(u)^{-1}$ is applied from the left (right).
One can generalize  {\em ordered products} \eqref{csym2} in an obvious way. 
For example,  
\beq\label{reversed}
K^{(n,m)}=\prod_{i=1,\ldots,n}^{\longleftarrow} \prod_{j=n+1,\ldots,n+m}^{\longleftarrow} \wvr{R}_{ij}(u_i+v_{j-n}-hc/2-hN)^{-1}
\eeq
satisfies
\beq\label{csym3}
{}^{rl\hspace{-3pt}}\left(K^{(n,m)}\right)\cdot\wndr{\wvr{R}}_{nm}^{12}(u-hc/2|v)=1,
\eeq
where superscript $rl$ in \eqref{csym3} indicates that the  tensor factors of  $K^{(n,m)}$ corresponding to the first index $i=1,\ldots,n$ in \eqref{reversed}  are applied from the right   in reversed order, while the tensor factors corresponding to the second index  $j=n+1,\ldots,n+m$ in \eqref{reversed}  are applied from the left.

\begin{ex}
Set 
$K_{ij}=\wvr{R}_{ij}(u_i+v_{j-n}-hc/2-hN)^{-1}$  and $S_{ij}=R_{ij}(u_i+v_{j-n}-hc/2)$.
We  briefly explain how to verify \eqref{csym3}
for $n=m=2$; the general case can be proved analogously. 
First, due to \eqref{rmatrixces} and \eqref{reversed}, on $(\ndo\CC^N)^{\ot 2}\ot (\ndo\CC^N)^{\ot 2}$ we have 
$$\wndr{\wvr{R}}_{22}^{12}(u-hc/2|v)=S_{13}S_{14}S_{23}S_{24}\fand K^{(2,2)}=K_{24}K_{23}K_{14}K_{13}. $$
The element ${}^{rl\hspace{-3pt}}\left(K^{(2,2)}\right)\cdot\wndr{\wvr{R}}_{22}^{12}(u-hc/2|v)$ can be written as
\beq\label{writtenas}{}^{rl\hspace{-3pt}}\left(K_{23}\right)\cdot
\left({}^{rl\hspace{-3pt}}\left(K_{24}\right)\cdot
\left({}^{rl\hspace{-3pt}}\left(K_{13}\right)\cdot
\left( {}^{rl\hspace{-3pt}}\left(K_{14}\right)\cdot
\left(
S_{13}S_{14}S_{23}S_{24}
\right)\right)
\right)
\right).
\eeq
By the first equality in \eqref{csym2} we have
$$ {}^{rl\hspace{-3pt}}\left(K_{14}\right)\cdot
\left(
S_{13}S_{14}S_{23}S_{24}
\right)=
S_{13}\left({}^{rl\hspace{-3pt}}\left(K_{14}\right)\cdot S_{14} \right)S_{23}S_{24}
=S_{13}S_{23}S_{24}.
$$
Next, as before, by the first equality in \eqref{csym2} we have
$${}^{rl\hspace{-3pt}}\left(K_{13}\right)\cdot
\left(
S_{13}S_{23}S_{24}
\right)
=
\left({}^{rl\hspace{-3pt}}\left(K_{13}\right)\cdot
S_{13}\right) S_{23}S_{24}
=S_{23}S_{24}.
$$
Hence, \eqref{writtenas} is equal to 
${}^{rl\hspace{-3pt}}\left(K_{23}\right)\cdot
\left({}^{rl\hspace{-3pt}}\left(K_{24}\right)\cdot
\left(
 S_{23}S_{24}
\right)
\right).
$
By repeating the same arguments two more times, we finally obtain ${}^{rl\hspace{-3pt}}\left(K^{(2,2)}\right)\cdot\wndr{\wvr{R}}_{22}^{12}(u-hc/2|v)=1$, as required.
\end{ex}

Observe that, due to \eqref{csym2}, the element 
\beq\label{haen}
L^{(n,m)}=\prod_{i=n+1,\ldots,n+m-1}^{\longleftarrow}\prod_{j=i+1,\ldots,n+m}^{\longleftarrow} \wvr{R}_{ij}(v_{i-n}+v_{j-n}-hN)^{-1}
\eeq
satisfies
\beq\label{csym4}
{}^{rl\hspace{-3pt}}\left(L^{(n,m)}\right)\cdot \wndr{B}_{m}^{+23}(v)=B_{n+1}^+(v_1)B_{n+2}^+(v_2)\ldots B_{n+m}^+(v_m),
\eeq
where, as before, superscript $rl$ in \eqref{csym4} indicates that the  tensor factors of $L^{(n,m)}$ corresponding to the first index $i=n+1,\ldots,n+m-1$ in \eqref{haen}  are applied from the right  in reversed order, while the tensor factors  corresponding to the second index  $j=i+1,\ldots,n+m$ in \eqref{haen} are applied from the left.
Relation \eqref{rbrb3}, together with \eqref{csym3} and \eqref{csym4}, implies the following equality  on
$(\ndo\CC^N)^{\ot n}\ot(\ndo\CC^N)^{\ot m}\ot \Ac_c(\gl_N)$:
\begin{align}
&\wndr{B}_{n}^{13}(u)\, B_{n+1}^+(v_1)B_{n+2}^+(v_2)\ldots B_{n+m}^+(v_m)={}^{rl\hspace{-3pt}}\left(L^{(n,m)}\right)\cdot\bigg({}^{rl\hspace{-3pt}}\left(K^{(n,m)}\right)\bigg.
\label{csym5}\\
&\qquad\bigg.\cdot\bigg( \wvr{R}_{nm}^{12}(u+3hc/2|v)^{-1}\wndr{B}_{m}^{+23}(v)\wndr{\wvr{R}}_{nm}^{12}(u+3hc/2|v)\wndr{B}_{n}^{13}(u)\wvr{R}_{nm}^{12}(u-hc/2|v)\bigg)\bigg).\non
\end{align}

Denote by $\vac$ the image of the unit $1\in\DY_c(\gl_N)$ in the quotient \eqref{imageofone}.
Let $\Wc_c '(\gl_N)$ be the $\B^+(\gl_N)$-submodule of $\Vc_c(\gl_N)$ generated by $\vac$.
Introduce the {\em vacuum module} $\Wc_c (\gl_N)$ over the algebra $\Ac_c(\gl_N)$ as the $h$-adic completion of $\Wc_c '(\gl_N)$.
Observe that $\Wc_c(\gl_N)$ is closed under the action of $\B_c(\gl_N)$, so it possesses a structure of an $\Ac_c(\gl_N)$-module.  Indeed, by applying \eqref{csym5} with $n=1$ on $\vac $
and using
\beq\label{invs}
B(u)\vac=T(u+hc)\ts G\ts T(-u)^{-1}\vac=T(u+hc)\ts G\vac=G\vac,
\eeq
we obtain 
\begin{align}
&B_{1}(u_1)\,B_{2}^+(v_1)B_{3}^+(v_2)\ldots B_{m+1}^+(v_m)\vac={}^{rl\hspace{-3pt}}\left(L^{(1,m)}\right)\cdot\bigg({}^{rl\hspace{-3pt}}\left(K^{(1,m)}\right)\bigg.\label{actionbu}\\
&\qquad\bigg.\cdot\bigg( \wvr{R}_{1m}^{12}(u_1 +3hc/2|v)^{-1}\wndr{B}_{m}^{+23}(v)\wndr{\wvr{R}}_{1m}^{12}(u_1 +3hc/2|v) \hspace{1pt} G_1\hspace{1pt}  \wvr{R}_{1m}^{12}(u_1 -hc/2|v)\vac\bigg)\bigg),\non
\end{align}
so it remains to observe that all coefficients of the right hand side in \eqref{actionbu} belong to $\Wc_c (\gl_N)$.

By the Poincar\'e--Birkhoff--Witt theorem for the double Yangian, see \cite[Theorem 2.2]{JKMY},
 the $\CC[[h]]$-modules $\Wc_c '(\gl_N)$ and  $\B^+(\gl_N)$ are isomorphic. Hence, in particular, the completion $\Wc_c (\gl_N)$ is topologically free, i.e. separated, torsion-free and $h$-adically complete.

\section{Quasi modules for $h$-adic nonlocal vertex algebras}\label{sec2}
In this section, we  study  $h$-adic nonlocal vertex algebras and their quasi modules, as defined by Li in \cite{Li}, and we establish some technical results on their center, which will be useful in Section \ref{sec3}. Next, we recall  Etingof--Kazhdan's definition \cite{EK} of  quantum VOA structure on the vacuum module $\Vc_{c}(\gl_N)$, $c\in\CC$. Finally, we construct   quasi modules  for the quantum VOA $\Vc_{2c}(\gl_N)$ on the $\CC[[h]]$-module $\Wc_{c}(\gl_N)$.

\subsection{Quasi modules}\label{thissection}
Let us recall the notion of quasi module for $h$-adic nonlocal vertex algebra; see \cite{Li}.  The tensor products in the next two definitions are   $h$-adically completed.
\begin{defn}
An {\em $h$-adic nonlocal vertex algebra}  is a triple $(V,Y,\vac)$, where $V$ is a topologically free $\mathbb{C}[[h]]$-module,
$\vac$ is a distinguished element of $V$ ({\em vacuum vector}) and 
\begin{align*}
Y \colon V\ot V&\to V((z))[[h]]\\
v\ot w&\mapsto Y(z)(v\ot w)=Y(v,z)w=\sum_{r\in\mathbb{Z}} v_r w z^{-r-1}
\end{align*}
  is a $\CC[[h]]$-module map
which satisfies the  {\em   weak associativity} property:
For any integer $n\geqslant 0$ and elements $u,v,w\in V$ there exists an   integer $r\geqslant 0$ such that
\beq\label{associativity}
(z_0 +z_2)^r  Y(v,z_0 +z_2)Y(w,z_2)  u - (z_0 +z_2)^r  Y\big(Y(v,z_0)w,z_2\big)  u\,\,
\in\,\, h^n V[[z_0^{\pm 1},z_2^{\pm 1}]];
\eeq
and the following conditions hold:
\begin{align*}
Y(v,z)\vac\in V[[z]],\quad \lim_{z\to 0}Y(v,z)\vac=v\fand Y(\vac,z)v=v\qquad\text{for any }v\in V.
\end{align*}
\end{defn}

\begin{defn}
Let $(V,Y,\vac)$  be an $h$-adic nonlocal vertex algebra. {\em Quasi $V$-module} is a pair $(W,Y_W)$, where $W$ is a topologically free $\CC[[h]]$-module
and 
\begin{align*}
Y_W(z)\colon V\ot W&\to W((z))[[h]]\\
v\ot w&\mapsto Y_W(z)(v\ot w)=Y_W(v,z)w=\sum_{r\in\mathbb{Z}} v_r w z^{-r-1}
\end{align*}
is a $\CC[[h]]$-module map
which satisfies the following:
For any integer $n\geqslant 0$ and elements $u,v\in V$, $w\in W$
there exists a nonzero polynomial $p(x_1,x_2)$ in $\CC[x_1,x_2]$
such that
\begin{align}
&p(z_0 +z_2,z_2)  Y_W(u,z_0 +z_2)Y_W(v,z_2)  w \non\\
&\qquad\qquad -   p(z_0 +z_2,z_2)Y_W\big(Y(u,z_0)v,z_2\big)w\,\,
\in\,\, h^n W[[z_0^{\pm 1},z_2^{\pm 1}]];\label{mod2}
\end{align}
and for any $w\in W$ we have $Y_W(\vac,z)w=w$. 
\end{defn}

Let $W$ be a  $\CC[[h]]$-module. For any  $a,b\in W[[z_0^{\pm 1},z_1^{\pm 1},\ldots]]$ and $n\geqslant 0$ we will write 
$$a\hsym b\qquad \text{if} \qquad a-b\,\in\, h^n W[[z_0^{\pm 1},z_1^{\pm 1},\ldots]].$$
\begin{lem}\label{lemaa}
Let $V$ be an $h$-adic nonlocal vertex algebra and let $W$ be a quasi $V$-module. 
Suppose that the elements $a,b\in V$ and $w_1,w_2\in W$ satisfy
\beq\label{ykom4}
[Y_{W}(a,z_1),Y_{W}(b,z_2)]w_i=0\quad\text{for } i=1,2.
\eeq
Then, for any  integers $p,t,n$,   $n\geqslant 0$, there exist scalars $\alpha_{r,s}\in\CC$, which do not depend on $i=1,2$, such that
\beq\label{ykom5}
(a_{p} b)_t w_i \,\hsym\, \sum_{r,s\in\mathbb{Z}} \alpha_{r,s} a_r b_s w_i\quad\text{for } i=1,2.
\eeq
\end{lem}

\begin{prf}
Fix integers $p,t,n$,   $n\geqslant 0$.
By \eqref{mod2},  there exist nonzero polynomials $p_i(x_1,x_2)$ in $\CC[x_1,x_2]$, where $i=1,2$, such that
\begin{equation}\label{pkom5}
p_i (z_0 +z_2,z_2) Y_{W}(a,z_0 +z_2)Y_{W}(b,z_2)  w_i \,\hsym\,   p_i(z_0 +z_2,z_2)Y_{W}\big(Y(a,z_0)b,z_2\big)w_i.
\end{equation}
Consider the left hand side in \eqref{pkom5}. Due to \eqref{ykom4},  
there exist an integer $l\geqslant 0$ such that
\beq\label{el}
(z_0 +z_2)^l z_2^l Y_{W}(a,z_0 +z_2)Y_{W}(b,z_2)  w_i =X_i(z_0,z_2) + h^n Z_i(z_0,z_2),\quad i=1,2,
\eeq
for some $X_i(z_0,z_2)\in W[[z_0,z_2]]$ and $Z_i(z_0,z_2)\in W((z_0))((z_2))[[h]]$.
Indeed, we can set $l=\max\left\{l_1,l_2,k_1,k_2 \right\}$, where $l_i$ and $k_i$ are chosen so that the expression 
$$z_1^{l_i}z_2^{k_i}Y_{W}(a,z_1)Y_{W}(b,z_2)w_i=z_1^{l_i}z_2^{k_i}Y_{W}(b,z_2)Y_{W}(a,z_1)w_i,\quad i=1,2,$$
possesses only nonnegative powers of the variables $z_1,z_2$ modulo $h^n$.
Equality \eqref{el} implies that   there exist scalars $\beta_{r,s}\in\CC$, which do not depend on $i=1,2$, such that the coefficient of $z_0^{-p-1}z_2^{-t-1}$ in $X_i(z_0,z_2)$ is equal to
\beq\label{ivana}
\sum_{r,s}\beta_{r,s} a_r b_s w_i \mod h^n\qquad \text{for }i=1,2.
\eeq

By combining relations \eqref{pkom5} and \eqref{el} we obtain
\beq\label{el2}
p_i(z_0 +z_2,z_2) X_i(z_0,z_2)\, \hsym\,  p_i(z_0 +z_2,z_2) (z_0 +z_2)^l z_2^l Y_{W}\big(Y(a,z_0)b,z_2\big)w_i,\quad i=1,2.
\eeq
The left hand side in \eqref{el2}, as well as $X_i(z_0,z_2) $, possesses only nonnegative powers of the variables  $z_0$ and $z_2$, while the right hand side in \eqref{el2}, as well as the expression 
$Y_{W}\big(Y(a,z_0)b,z_2\big)w_i$, belongs to $W ((z_2))((z_0))[[h]]$.
Hence, we can multiply \eqref{el2} by the inverse of the polynomial $p(z_2+z_0,z_2)$ in $\CC((z_2))((z_0))$, thus getting
\beq\label{el3}
X_i(z_0,z_2) \,\hsym\,   (z_0 +z_2)^l z_2^l Y_{W}\big(Y(a,z_0)b,z_2\big)w_i,\quad i=1,2.
\eeq
Next, we multiply \eqref{el3} by the inverse of the polynomial $(z_2 +z_0)^l z_2^l$ in $\CC((z_2))((z_0))$, which gives us
\beq\label{el4}
\left((z_2 +z_0)^l z_2^l\right)^{-1}\cdot X_i(z_0,z_2) \,\hsym\,    Y_{W}\big(Y(a,z_0)b,z_2\big)w_i,\quad i=1,2.
\eeq
In particular,  the coefficients of $z_0^{-p-1}z_2^{-t-1}$ in \eqref{el4} coincide modulo $h^n$.
Recall \eqref{ivana}. Clearly,  there exist scalars $\alpha_{r,s}\in\CC$, which do not depend on $i=1,2$, such that the coefficient of
$z_0^{-p-1} z_2^{-t-1}$ on the left hand side in \eqref{el4}   equals
$$\sum_{r,s\in\ZZ}\alpha_{r,s}a_rb_s w_i\mod h^n \qquad \text{for }i=1,2.$$
Since the coefficient of $z_0^{-p-1}$ on the right hand side in \eqref{el4}  equals $Y_{W}(a_{p} b,z_2\big)w_i$, by taking the coefficient of $z_2^{-t-1}$ we obtain
$$\sum_{r,s\in\ZZ}\alpha_{r,s}a_rb_s w_i\,\hsym\, (a_{p} b)_t w_i ,\quad i=1,2,$$
as required.
\end{prf}

As with quantum VOAs in \cite{JKMY}, we can introduce the center of an  $h$-adic nonlocal vertex algebra $V$ in analogy with vertex algebra theory; see, e.g., \cite[Chapter 2]{F}.  Define the {\em center}  of $V$ as the $\CC[[h]]$-submodule 
$$\z(V)=\left\{v\in V\,:\, w_r v=0\text{ for all }w\in V\text{ and }r\geqslant 0   \right\}.$$
It is worth noting that, in contrast with the vertex algebra theory, the center of an 
$h$-adic nonlocal vertex algebra does not need to be commutative; see \cite[Proposition 4.2]{JKMY}.

\begin{pro}\label{protech}
Let $V$ be an $h$-adic nonlocal vertex algebra and let $W$ be a quasi $V$-module. Suppose that the center $\z(V)$ is a commutative associative algebra, with respect to the product $a\cdot b\coloneqq a_{-1}b$ for $a,b\in \z(V)$. Furthermore, assume that the algebra $\z(V)$ is topologically generated, with respect to the $h$-adic topology, by some family $\Phi\subseteq \z(V)$.
\begin{enumerate}[(a)]
\item\label{aaaqqq}
If
$[Y_W(a, z_1),Y_W(b,z_2) ]=0$  for all $a,b\in \Phi$,
then
\beq\label{1com}[Y_W(a, z_1),Y_W(b,z_2) ]=0\qquad\text{for all }a,b\in\z(V).\eeq
\item\label{bbbqqq}
If $\psi\colon W\to W$ is a $\CC[[h]]$-module map satisfying $[Y_W(a, z),\psi ]=0$  for all $a\in \Phi$, then
\beq\label{2com}
[Y_W(a, z),\psi ]=0\qquad\text{for all }a \in\z(V).
\eeq
\end{enumerate}
\end{pro}

\begin{prf}
Let $a,b,c$ be elements of the center $\z(V)$ such that the pairs $(a,b)$, $(b,c)$ and $(a,c)$ satisfy \eqref{1com}. In order to prove   \eqref{aaaqqq}, it is sufficient to verify that the pair $(a\cdot b,c)=(a_{-1}b,c)$ satisfies \eqref{1com}.
 Fix  $w\in W$ and integers $p,t,n$,   $n\geqslant 0$. 
Due to our assumption, the pair $(a,b)$ satisfies \eqref{1com}, so  Lemma \ref{lemaa} implies that there exist scalars $\alpha_{r,s}\in\CC$ such that
\beq\label{ppp0}
(a\cdot b)_t c_p w\,\hsym\,\sum_{r,s\in\ZZ} \alpha_{r,s} a_r b_s c_p w 
\Fand (a\cdot b)_t w\,\hsym\,\sum_{r,s\in\ZZ} \alpha_{r,s} a_r b_s  w.
\eeq
Since the pairs $(a,c)$ and $(b,c)$ satisfy \eqref{1com}, we have
$[b_s,c_p]=[a_r,c_p]=0$ on $W$, 
 so by  relations in \eqref{ppp0} we have
$$
(a\cdot b)_t c_pw\hsym\sum_{r,s\in\ZZ} \alpha_{r,s} a_r b_s c_p w= \sum_{r,s\in\ZZ} \alpha_{r,s} a_r  c_p b_s w= \sum_{r,s\in\ZZ} \alpha_{r,s}   c_p a_r b_s w\hsym c_p(a\cdot b)_t w.
$$
Hence we proved
$
(a\cdot b)_t c_p w\,\hsym\,  c_p(a\cdot b)_t w.
$
Since $n$ was arbitrary, we conclude that
$$
(a\cdot b)_t c_pw=  c_p(a\cdot b)_t w.
$$
Finally, since integers $p,t$ and element $w\in W$ were arbitrary, this gives us
$$
[Y_{W}(a\cdot b,z_1),Y_{W}(c,z_2)]=0,
$$
as required. 
Statement \eqref{bbbqqq} can be proved analogously.
\end{prf}

\subsection{Vacuum module \texorpdfstring{$\Vc_c (\gl_N)$}{Vc(glN)} as a quantum VOA}

Let $n$ and $m$ be positive integers. For the families of variables $u=(u_1,\ldots,u_n)$ and $v=(v_{1},\ldots, v_{m})$ and the variable $z$ consider
the functions  with values in
$(\ndo\CC^N)^{\ot n} \ot (\ndo\CC^N)^{\ot m}$
\begin{align}
& \wvr{R}_{nm}^{12}(u|v|z)=\prod_{i=1,\ldots,n}^{\longrightarrow} \prod_{j=n+1,\ldots,n+m}^{\longleftarrow}\wvr{R}_{ij}(z+u_i-v_{j-n}),\label{bnm1}\\
&\wndr{\wvr{R}}_{nm}^{12}(u|v|z)=\prod_{i=1,\ldots,n}^{\longrightarrow} \prod_{j=n+1,\ldots,n+m}^{\longrightarrow}\wvr{R}_{ij}(z+u_i+v_{j-n}).\label{bnm2}
\end{align}
The functions   $R_{nm}^{12}(u|v|z)$ and $\wndr{R}_{nm}^{12}(u|v|z)$  corresponding to  $R$-matrix \eqref{yang} can be defined analogously.
In \eqref{bnm1}--\eqref{bnm2}, as well as in the rest of the paper, we use the common expansion convention: expressions of the form $(a_1 z_1+\ldots+ a_n z_n )^k$, where $a_i \in\CC$, $a_i\neq 0$ and $k<0$, are expanded in negative powers of the variable appearing on the left, e.g., 
$$(z_1 -z_2)^{-1}=\sum_{l\geqslant 0}\frac{z_2^l}{z_1^{l+1}}\in\CC[z^{-1}_1][[z_2]] \fand (-z_2+z_1)^{-1}=-\sum_{l\geqslant 0}\frac{z_1^l}{z_2^{l+1}}\in\CC[z^{-1}_2][[z_1]].$$
In particular, \eqref{bnm1}--\eqref{bnm2} contain only nonnegative powers of the variables $u_i$ and $v_{j-n}$.

Define the following operators on $(\ndo\CC^N)^{\ot n} \ot\Vc_c(\gl_N)$:
\begin{align*}
T_n^+(u|z)=T_1^+(z+u_1)\ldots T_n^+(z+u_n)\fand T_n(u|z)=T_1(z+u_1)\ldots T_n(z+u_n).
\end{align*}
Using  \eqref{RTT2}--\eqref{RTT3},  one can easily verify the following equations for the operators on $(\ndo\CC^N)^{\ot n} \ot (\ndo\CC^N)^{\ot m}\ot \Vc_c(\gl_N)$, originally given in \cite{EK}, which employ the superscript notation introduced prior to Proposition \ref{rbrb_pro}.
\begin{align}
&R_{nm}^{12}(u|v|z-w)T_n^{+13}(u|z)T_m^{+23}(v|w)
=T_m^{+23}(v|w)T_n^{+13}(u|z)R_{nm}^{12}(u|v|z-w),\label{rtt1}\\ 
&R_{nm}^{12}(u|v|z-w)T_n^{13}(u|z)T_m^{23}(v|w)
=T_m^{23}(v|w)T_n^{13}(u|z)R_{nm}^{12}(u|v|z-w),\label{rtt2}\\ 
&\overline{R}_{nm}^{\ts 12}(u|v|z-w+h\tss c/2)T_n^{13}(u|z)T_m^{+23}(v|w)\nonumber\\
&\qquad\qquad\qquad\qquad\qquad=T_m^{+23}(v|w)T_n^{13}(u|z)\overline{R}_{nm}^{\ts 12}(u|v|z-w-h\tss c/2).
\label{rtt3}
\end{align}

The next theorem, which  is due to Etingof and Kazhdan \cite{EK}, introduces the structure of   quantum VOA on the vacuum module $\Vc_c(\gl_N)$.
Roughly speaking, quantum VOA, as defined in \cite{EK}, is an $h$-adic nonlocal vertex algebra $(V,Y,\vac)$  equipped with the
$\mathbb{C}[[h]]$-module map $\mathcal{S}(z)\colon V\ot V\to V\ot V\otimes\mathbb{C}((z))$ (with the tensor products being   $h$-adically completed) satisfying the
{\em $\mathcal{S}$-locality}:
For any integer $n\geqslant 0$ and elements $v,w\in V$ there exists
an integer 
$k\geqslant 0$ such that for any $u\in V$
\begin{align}
&(z_1-z_2)^{k} Y(z_1)\big(1\otimes Y(z_2)\big)\big(\mathcal{S}(z_1 -z_2)(v\otimes w)\otimes u\big)\nonumber\\
&\qquad-(z_1-z_2)^{k} Y(z_2)\big(1\otimes Y(z_1)\big)(w\otimes v\otimes u)
\,\in\, h^n V[[z_1^{\pm 1},z_2^{\pm 1}]];\label{locality}
\end{align}
and several other properties. 
In this paper, we only use $\mathcal{S}$-locality \eqref{locality} and the underlying structure of an $h$-adic nonlocal vertex algebra on $\Vc_c(\gl_N)$, so we omit the original definition of quantum VOA.  

\begin{thm}\label{EK:qva}
For any $c\in \CC$
there exists a unique  structure of quantum VOA
on  $\Vc_c(\gl_N)$ such that the vacuum vector is
$\vac\in \Vc_c(\gl_N)$, the vertex operator map is given by
\beq\label{qva1}
Y\big(T_n^{+}(u|0)\vac,z\big)=T_n^{+}(u|z)\ts T_n (u|z+h\tss c/2)^{-1}
\eeq
and the map $\mathcal{S}(z)$ is defined by the relation
\begin{align}
\mathcal{S}^{34}(z)&\Big(\overline{R}_{nm}^{  12}(u|v|z)^{-1}  T_{m}^{+24}(v|0) 
\overline{R}_{nm}^{  12}(u|v|z-h  c)  T_{n}^{+13}(u|0)(\vac\otimes \vac) \Big)\nonumber\\
 =&T_{n}^{+13}(u|0)  \overline{R}_{nm}^{  12}(u|v|z+h  c)^{-1} 
T_{m}^{+24}(v|0)  \overline{R}_{nm}^{  12}(u|v|z)(\vac\otimes \vac)\label{qva4}
\end{align}
for operators on
$
(\ndo\mathbb{C}^{N})^{\otimes n} \otimes
(\ndo\mathbb{C}^{N})^{\otimes m}\otimes \Vc_c(\gl_N) \otimes \Vc_c(\gl_N) $.
\end{thm}

\subsection{Vacuum module \texorpdfstring{$\Wc_c (\gl_N)$}{Wc(glN)} as a quasi \texorpdfstring{$\Vc_{2c} (\gl_N)$}{V2c(glN)}-module}

Consider the operators on $(\ndo\CC^N)^{\ot n}\ot \Wc_c(\gl_N)$ given by
\begin{align*}
&\wndr{B}_{n}^{+}(u|z)=\prod_{i=1,\ldots,n}^{\longrightarrow} \left(B_{i}^+ (z+u_i) \wvr{R}_{i\, i+1}(2z+u_i+u_{i+1})\ldots \wvr{R}_{in}(2z+u_i+u_n)  \right)\fand\\
&\wndr{B}_{n}(u|z)=\prod_{i=1,\ldots,n}^{\longrightarrow} \left(B_{i} (z+u_i) \wvr{R}_{i\, i+1}(2z+u_i+u_{i+1}+hc)\ldots \wvr{R}_{in}(2z+u_i+u_n+hc)   \right).
\end{align*}
The next proposition can be proved by using Proposition \ref{rbrb_pro}.

\begin{pro}\label{rbrb_pro2}
Let $n$ and $m$ be positive integers. The following equalities hold for the operators on
$(\ndo\CC^N)^{\ot n}\ot(\ndo\CC^N)^{\ot m}\ot \Wc_c(\gl_N)$:
\begin{align}
&R_{nm}^{12}(u|v|z-w)\wndr{B}_{n}^{+13}(u|z)\wndr{R}_{nm}^{12}(u|v|z+w)\wndr{B}_{m}^{+23}(v|w)\non\\
&\quad=\wndr{B}_{m}^{+23}(v|w)\wndr{R}_{nm}^{12}(u|v|z+w)\wndr{B}_{n}^{+13}(u|z)R_{nm}^{12}(u|v|z-w),\label{rbrb12}\\
&R_{nm}^{12}(u|v|z-w)\wndr{B}_{n}^{13}(u|z)\wndr{R}_{nm}^{12}(u|v|z+w+hc)\wndr{B}_{m}^{23}(v|w)\non\\
&\quad=\wndr{B}_{m}^{23}(v|w)\wndr{R}_{nm}^{12}(u|v|z+w+hc)\wndr{B}_{n}^{13}(u|z)R_{nm}^{12}(u|v|z-w),\label{rbrb22}\\
&\wvr{R}_{nm}^{12}(u|v|z-w+3hc/2)\wndr{B}_{n}^{13}(u|z)\wndr{\wvr{R}}_{nm}^{12}(u|v|z+w-hc/2)\wndr{B}_{m}^{+23}(v|w)\non\\
&\quad=\wndr{B}_{m}^{+23}(v|w)\wndr{\wvr{R}}_{nm}^{12}(u|v|z+w+3hc/2)\wndr{B}_{n}^{13}(u|z)\wvr{R}_{nm}^{12}(u|v|z-w-hc/2).\label{rbrb32}
\end{align}
\end{pro}

The following theorem is our main result.

\begin{thm}\label{zmaina}
For any $c\in\CC$ there exists a unique structure of quasi $\Vc_{2c}(\gl_N)$-module  on the vacuum module $\Wc_c(\gl_N)$ such that
\beq\label{Y}
Y_{\Wc_c(\gl_N)}(T_n^+(u|0)\vac,z)=
\wndr{B}_{n}^{+}(u|z)\wndr{B}_{n}(u|z+hc/2)^{-1}.
\eeq
\end{thm}

\begin{prf}
Set $\Wc_c=\Wc_c(\gl_N)$.
We first prove that    map \eqref{Y} is well-defined. It is sufficient to verify that   $a\mapsto Y_{\Wc_c}(a,z)$  maps the ideal of   relations \eqref{RTT1} to itself since, due to  Poincar\'e--Birkhoff--Witt theorem for the double Yangian \cite[Theorem 2.2]{JKMY}, $\Y^+ (\gl_N)$ is isomorphic to the algebra generated by the elements $t_{ij}^{(-r)}$, where $r=1,2\ldots$ and $i,j=1,\ldots ,N$, subject to \eqref{RTT1}.  Set $\widehat{R}_{k\tss k+1}=R_{k\tss k+1}(u_k -u_{k+1})$, where $1\leqslant k<n$. Relation \eqref{RTT1} implies
\begin{align}
&\widehat{R}_{k\, k+1}T_n^+(u|0)\vac\non\\
=&T_{1}^+(u_{1})\ldots T_{k-1}^+(u_{k-1})T_{k+1}^+(u_{k+1})T_{k}^+(u_{k})T_{k+2}^+(u_{k+2})\ldots T_{n}^+(u_{n})\vac\widehat{R}_{k\tss k+1}.\label{prf0}
\end{align}
Set
$\wtld{R}_{ij}^z=\wvr{R}_{ij}(2z+u_i+u_j)$ and $\wtld{R}_{ij}^{z+hc}=\wvr{R}_{ij}(2z+u_i+u_j+2hc)$.
Due to Yang--Baxter equation \eqref{ybe} and unitarity \eqref{uni}, for any indices $1\leqslant j<k<k+1<l\leqslant n$ we have
\begin{align}
\widehat{R}_{k\tss k+1} \wtld{R}_{jk}^z\wtld{R}_{j\tss k+1}^z=\wtld{R}_{j\tss k+1}^z  \wtld{R}_{jk}^z \widehat{R}_{k\tss k+1}\fand
\widehat{R}_{k\tss k+1} \wtld{R}_{kl}^z\wtld{R}_{k+1\tss l}^z=\wtld{R}_{k+1\tss l}^z  \wtld{R}_{kl}^z \widehat{R}_{k\tss k+1}.\label{prf2}
\end{align}
Relation \eqref{rbrb12}, together with  \eqref{prf2}, implies
\beq\label{prf3}
\widehat{R}_{k\tss k+1}\wndr{B}_{n}^{+}(u|z)=\wndr{B}_{n, k\leftrightarrow k+1}^{+}(u|z)\widehat{R}_{k\tss k+1},\qquad\text{where}
\eeq
\begin{align*}
&\wndr{B}_{n, k\leftrightarrow k+1}^{+}(u|z) =\prod_{i=1,\ldots,k-1}^{\longrightarrow} \left(B_{i}^+ (z+u_i) \wtld{R}_{i\ts i+1}^z\ldots \wtld{R}_{i\tss k-1}^z \wtld{R}_{i\tss k+1}^z \wtld{R}_{i\tss k}^z \wtld{R}_{i\tss k+2}^z \ldots \wtld{R}_{i\tss n}^z    \right)\\
&\qquad\cdot\left(B_{k+1}^+ (z+u_{k+1}) \wtld{R}_{k+1\ts k}^z  \wtld{R}_{k+1\tss k+2}^z  \ldots \wtld{R}_{k+1\tss n}^z    \right) \cdot\left(B_{k}^+ (z+u_{k}) \wtld{R}_{k\ts k+2}^z    \ldots \wtld{R}_{k\tss n}^z    \right)\\
&\qquad\cdot
\prod_{i=k+2,\ldots,n}^{\longrightarrow} \left(B_{i}^+ (z+u_i) \wtld{R}_{i\ts i+1}^z\ldots \wtld{R}_{i\tss n}^z    \right).
\end{align*}

Next, due to Yang--Baxter equation \eqref{ybe} and unitarity \eqref{uni} we have 
\begin{align}
   \widehat{R}_{k\tss k+1}(\wtld{R}_{k+1\tss l}^{z+hc})^{-1} (\wtld{R}_{kl}^{z+hc})^{-1}&=(\wtld{R}_{kl}^{z+hc} )^{-1}(\wtld{R}_{k+1\tss l}^{z+hc})^{-1} \widehat{R}_{k\tss k+1},\label{prfy}\\
   \widehat{R}_{k\tss k+1}(\wtld{R}_{j\tss k+1}^{z+hc})^{-1}(\wtld{R}_{jk}^{z+hc})^{-1}&=(\wtld{R}_{jk}^{z+hc})^{-1}(\wtld{R}_{j\tss k+1}^{z+hc})^{-1} \widehat{R}_{k\tss k+1}\label{prfx}
\end{align}
for $1\leqslant j<k<k+1<l\leqslant n$. 
Relation \eqref{rbrb22}, together with \eqref{prfy}--\eqref{prfx}, implies
\beq\label{prf4}
\widehat{R}_{k\tss k+1}\wndr{B}_{n}(u|z+hc/2)^{-1}=\wndr{B}_{n, k\leftrightarrow k+1}(u|z+hc/2)^{-1}\widehat{R}_{k\tss k+1},\qquad\text{where}
\eeq
\begin{align*}
&\wndr{B}_{n, k\leftrightarrow k+1}(u|z+hc/2) =\prod_{i=1,\ldots,k-1}^{\longrightarrow} \bigg(B_{i} (z+u_i+hc/2) \wtld{R}_{i\ts i+1}^{z+hc}\ldots \wtld{R}_{i\tss k-1}^{z+hc}\bigg.\\
& \bigg.\cdot \wtld{R}_{i\tss k+1}^{z+hc} \wtld{R}_{i\tss k}^{z+hc} \wtld{R}_{i\tss k+2}^{z+hc} \ldots \wtld{R}_{i\tss n}^{z+hc}    \bigg)\cdot \left(B_{k+1} (z+u_{k+1}+hc/2) \wtld{R}_{k+1\ts k}^{z+hc}  \wtld{R}_{k+1\tss k+2}^{z+hc}  \ldots \wtld{R}_{k+1\tss n}^{z+hc}    \right)\\
&\cdot\left(B_{k} (z+u_{k}+hc/2) \wtld{R}_{k\ts k+2}^{z+hc}    \ldots \wtld{R}_{k\tss n}^{z+hc}    \right)\cdot
\prod_{i=k+2,\ldots,n}^{\longrightarrow} \left(B_{i} (z+u_i+hc/2) \wtld{R}_{i\ts i+1}^{z+hc}\ldots \wtld{R}_{i\tss n}^{z+hc}    \right).
\end{align*}

Finally, by applying the   map $a\mapsto Y_{\Wc_c}(a,z)$ on the left hand side of \eqref{prf0}, we obtain
$$\widehat{R}_{k\tss k+1}\wndr{B}_{n}^{+}(u|z)\wndr{B}_{n}(u|z+hc/2)^{-1}.$$
By \eqref{prf3} and \eqref{prf4}   this is equal to
\beq\label{prf6}
\wndr{B}_{n, k\leftrightarrow k+1}^{+}(u|z)\wndr{B}_{n, k\leftrightarrow k+1}(u|z+hc/2)^{-1}\widehat{R}_{k\tss k+1}.
\eeq
However, \eqref{prf6} coincides with the image of the right hand side in \eqref{prf0}, with respect to the map $a\mapsto Y_{\Wc_c}(a,z)$, so we conclude that
$Y_{\Wc_c}(z)$ is well-defined. 

It is clear that \eqref{Y} determines the   map $Y_{\Wc_c}(z)$  uniquely.
Our next goal is to show that the image of $Y_{\Wc_c}(z)$ belongs to $\Wc_c((z))[[h]]$.
Relation \eqref{rbrb32} implies
\begin{align}
&\wvr{R}_{nm}^{12}(u|v|z+2hc)\wndr{B}_{n}^{13}(u|z+hc/2)\wndr{\wvr{R}}_{nm}^{12}(u|v|z)\wndr{B}_{m}^{+23}(v)\non\\
&\qquad=\wndr{B}_{m}^{+23}(v)\wndr{\wvr{R}}_{nm}^{12}(u|v|z+2hc)\wndr{B}_{n}^{13}(u|z+hc/2)\wvr{R}_{nm}^{12}(u|v|z).\label{axiom1}
\end{align}
Observe that 
$$\wndr{B}_{n}^{13}(u|z+hc/2)^{-1}\vac=\wht{G}\vac\quad\text{for}\quad \wht{G}=G_1 \ldots G_n,$$
so, by   using \eqref{csym2} and \eqref{axiom1} and  arguing as in the proof of Equality  \eqref{csym5} we get
\begin{align}
&Y_{\Wc_c}(T_n^+(u|0)\vac,z)B_{n+1}^{+}(v_1)\ldots B_{n+m}^{+}(v_m)\vac\label{lbel}\\
=&\wndr{B}_{n}^{+13}(u|z)\wndr{B}_{n}^{13}(u|z+hc/2)^{-1} B_{n+1}^{+}(v_1)\ldots B_{n+m}^{+}(v_m)\vac\non\\
=&\wndr{B}_{n}^{+13}(u|z)\hspace{-2pt}\left(\hspace{-2pt} {}^{rl\hspace{-3pt}}\left(L\right)\cdot\hspace{-2pt}\left(\hspace{-2pt}{}^{rl\hspace{-3pt}}\left(K\right)\cdot\hspace{-2pt}\left(\wndr{\wvr{R}}_{nm}^{12}(u|v|z)\wndr{B}_{m}^{+23}(v)\wvr{R}_{nm}^{12}(u|v|z)^{-1}
\,\wht{G}\vac
\wndr{\wvr{R}}_{nm}^{12}(u|v|z+2hc)^{-1}\right)\hspace{-2pt}\right) \hspace{-2pt}\right),\non
\end{align}
where 
$K=  \wvr{R}_{nm}^{12}(u|v|z+2hc+hN)$
and $L=L^{(n,m)}$ is given by \eqref{haen}. 
Recall that the $R$-matrix $\wvr{R}(x)$ belongs to $(\ndo\CC^N) [x^{-1}][[h]]$. Therefore, the right hand side of \eqref{lbel} is a Taylor series in the variables $u_1,\ldots,u_n,v_1,\ldots,v_m$ and $h$ such that
the coefficient of each monomial 
$u_1^{a_1}\ldots u_n^{a_n}v_1^{b_1}\ldots v_m^{b_m}h^b$ 
possesses only finitely many negative powers of the variable $z$.  This implies that the image of   $Y_{\Wc_c}(z)$
belongs to $\Wc_c((z))[[h]]$.

The property $Y_{\Wc_c}(\vac,z)=1_{\Wc_c}$ is clear, so it remains to prove \eqref{mod2}. Consider  the second summand in \eqref{mod2}.
By applying  the vertex operator map $Y(z_0)$ for the quantum VOA $\Vc_{2c}(\gl_N)$, as defined in Theorem \ref{EK:qva}, on the series\footnote{It is possible (and perhaps more natural) to prove \eqref{mod2} by starting from   $T_n^{+13}(u|0) T_m^{+24}(v|0)(\vac\ot \vac)$ instead of \eqref{as1}. However, this requires the use of ordered products, as defined in Section \ref{explinedin}, thus making the calculations seemingly more complicated, even though the proof  remains analogous.}
\beq\label{as1}
T_n^{+13}(u|0)\wvr{R}^{12}_{nm}(u|v|z_0 +2hc)^{-1}T_m^{+24}(v|0)(\vac\ot \vac),
\eeq 
whose coefficients belong to $(\ndo \CC^N)^{\ot n}\ot (\ndo \CC^N)^{\ot m}\ot\Vc_{2c}(\gl_N)\ot\Vc_{2c}(\gl_N) $, we get
$$T_n^{+13}(u|z_0)T_n^{13}(u|z_0+hc)^{-1}\wvr{R}^{12}_{nm}(u|v|z_0 +2hc)^{-1}T_m^{+23}(v|0)\vac.$$
Due to \eqref{rtt3} at the level $2c$ and $T_n^{13}(u|z_0 +hc)^{-1}\vac=\vac$ this equals to
\beq\label{as2}
T_n^{+13}(u|z_0)T_m^{+23}(v|0)\vac \wvr{R}^{12}_{nm}(u|v|z_0 )^{-1},
\eeq
which is a series with coefficients in $(\ndo \CC^N)^{\ot n}\ot (\ndo \CC^N)^{\ot m}\ot\Vc_{2c}(\gl_N)$. Finally,  by applying the map $a\mapsto Y_{\Wc_c}(a,z_2)$ on \eqref{as2} we get
\beq\label{as3}
\wndr{B}_{n+m}^{+}(x|z_2)\wndr{B}_{n+m}(x|z_2+hc/2)^{-1}\wvr{R}^{12}_{nm}(u|v|z_0 )^{-1},
\eeq
where $x$ denotes the $n+m$ variables
$x= (z_0+u_1,\ldots,z_0+u_n,v_1,\ldots,v_m)$.

Let us consider the first summand in \eqref{mod2}. By applying
$Y_{\Wc_c}(z_0 +z_2)(1\ot Y_{\Wc_c}(z_2))$ on \eqref{as1} we obtain
\begin{align}
&\wndr{B}^{+13}_{n}(u|z_0+z_2)\wndr{B}_{n}^{13}(u|z_0+z_2+hc/2)^{-1}
\non\\
&\qquad\cdot\wvr{R}^{12}_{nm}(u|v|z_0 +2hc)^{-1}\wndr{B}^{+23}_{m}(v|z_2)\wndr{B}^{23}_{m}(v|z_2+hc/2)^{-1}.\label{as4}
\end{align}
Using \eqref{rbrb32} we can express   
$ \wndr{B}_{n}^{13}(u|z_0+z_2+hc/2)^{-1}
\wvr{R}^{12}_{nm}(u|v|z_0 +2hc)^{-1}
\wndr{B}^{+23}_{m}(v|z_2)$ 
  as
\begin{align*}
&\wndr{\wvr{R}}_{nm}^{12}(u|v|z_0 +2z_2)\wndr{B}_{m}^{+23}(v|z_2)\non\\
&\qquad\cdot\wvr{R}_{nm}^{12}(u|v|z_0)^{-1}\wndr{B}_{n}^{13}(u|z_0+z_2+hc/2)^{-1}\wndr{\wvr{R}}_{nm}^{12}(u|v|z_0+2z_2+2hc)^{-1},
\end{align*}
so  that \eqref{as4} is equal to
\begin{align}
&\wndr{B}^{+13}_{n}(u|z_0+z_2)\wndr{\wvr{R}}_{nm}^{12}(u|v|z_0 +2z_2)\wndr{B}_{m}^{+23}(v|z_2)\wvr{R}_{nm}^{12}(u|v|z_0)^{-1}\non\\
&\qquad\cdot\wndr{B}_{n}^{13}(u|z_0+z_2+hc/2)^{-1}\wndr{\wvr{R}}_{nm}^{12}(u|v|z_0+2z_2+2hc)^{-1}\wndr{B}_{m}^{23}(v|z_2+hc/2)^{-1}.\label{as5}
\end{align}
Finally, we   rewrite \eqref{as5} using \eqref{rbrb22}, thus getting
\begin{samepage}
\begin{align}
&\wndr{B}^{+13}_{n}(u|z_0+z_2)\wndr{\wvr{R}}_{nm}^{12}(u|v|z_0 +2z_2)\wndr{B}_{m}^{+23}(v|z_2)\wndr{B}_{m}^{23}(v|z_2+hc/2)^{-1}\non\\ 
&\qquad\cdot\wndr{\wvr{R}}_{nm}^{12}(u|v|z_0+2z_2+2hc)^{-1}\wndr{B}_{n}^{13}(u|z_0+z_2+hc/2)^{-1}\wvr{R}_{nm}^{12}(u|v|z_0)^{-1}.\label{as77}
\end{align}
\end{samepage}

Expressions \eqref{as3} and \eqref{as77} are not  equal, even though they do coincide  when viewed as Taylor series in the variables $u_1,\ldots ,u_n,v_1,\ldots ,v_m,h$ whose coefficients are rational functions in  $z_0,z_2$. Indeed, due to our expansion convention, the operators and $R$-matrices in \eqref{as3}, whose arguments contain both the variables $z_0$ and $z_2$, should be expanded in nonnegative powers of $z_0$, while the same operators and $R$-matrices in \eqref{as77} should be expanded in nonnegative powers of $z_2$. Fix an  integer $k\geqslant 0$ and an element $w\in\Wc_c $. Apply both  \eqref{as3} and \eqref{as77} on $w$ and denote the resulting expressions by $P(u,v,z_0,z_2)$ and $S(u,v,z_0,z_2)$ respectively. Then, for any choice of  integers
$a_1,\ldots, a_n\geqslant 0$ and $b_1,\ldots,b_m\geqslant 0$ there exist an  integer $r\geqslant 0$ such that the coefficients of
$u^{a_1}_1\ldots u^{a_n}_n v^{b_1}_1\ldots v^{b_m}_m$ in
$$(z_0+z_2)^r (z_0+2z_2)^r P(u,v,z_0,z_2)\fand (z_0+z_2)^r (z_0+2z_2)^r S(u,v,z_0,z_2)$$
coincide modulo $h^k$, which implies \eqref{mod2}.
\end{prf}

The   map $Y_{\Wc_c(\gl_N)}(z)$ satisfies the following "twisted" $\Sc$-locality property; cf. \cite{Li0,Li2}.
\begin{pro}\label{proploc}
For any  $u,v\in \Vc_{2c}(\gl_N)$ and   integer $k  \geqslant 0 $ there exists an integer
$r \geqslant 0 $ such that for any $w\in \Wc_c(\gl_N)$
\begin{align}
&(z_1^2-z^2_2)^r\, Y_{\Wc_c(\gl_N)}(z_1)\big(1\otimes Y_{\Wc_c(\gl_N)}(z_2)\big)\big(\mathcal{S}(z_1 -z_2)(u\otimes v)\otimes w\big)
\label{mod3}\\
&\qquad-(z^2_1-z^2_2)^r\, Y_{\Wc_c(\gl_N)}(z_2)\big(1\otimes Y_{\Wc_c(\gl_N)}(z_1)\big)(v\otimes u\otimes w)
\,\in\, h^k {\Wc_c}(\gl_N)[[z_1^{\pm 1},z_2^{\pm 1}]]. \non
\end{align}
\end{pro}

\begin{prf}
Set $\Wc_c=\Wc_c(\gl_N)$. 
Consider the first summand in \eqref{mod3} and
set
$z=z_1-z_2$. Notice that the variable $z_1$ appears on the left in $z=z_1-z_2$, so  the negative powers of $z$ should be expanded in negative powers of $z_1$.
By applying $\Sc (z)$ at the level $2c$, as defined in  \eqref{qva4}, on the last two tensor factors 
of the expression
\beq\label{loc1}
\wvr{R}_{nm}^{12}(u|v|z)^{-1}T_m^{+24}(v|0)\wvr{R}_{nm}^{12}(u|v|z-2hc)T_n^{+13}(u|0)(\vac\ot\vac),
\eeq
whose coefficients belong to $(\ndo \CC^{N})^{\ot (n+m)}\ot\Vc_{2c}(\gl_N)^{\ot 2}$,
we get
\beq\label{loc2}
T_n^{+13}(u|0)\wvr{R}_{nm}^{12}(u|v|z+2hc)^{-1}T_m^{+24}(v|0)\wvr{R}_{nm}^{12}(u|v|z)(\vac\ot\vac).
\eeq
Next, we apply $Y_{\Wc_c}(z_1)(1\ot Y_{\Wc_c}(z_2))$ on \eqref{loc2}, thus getting 
\begin{align}
&\wndr{B}_n^{+13}(u|z_1)\wndr{B}_n^{13}(u|z_1+hc/2)^{-1}\wvr{R}_{nm}^{12}(u|v|z+2hc)^{-1}\non\\
&\qquad\qquad\cdot\wndr{B}_m^{+23}(v|z_2)\wndr{B}_m^{23}(v|z_2+hc/2)^{-1}\wvr{R}_{nm}^{12}(u|v|z).\label{loc3}
\end{align}
We may now proceed as in  calculation \eqref{as4}--\eqref{as77} and prove that \eqref{loc3} equals
\begin{align}
&\wndr{B}_n^{+13}(u|z_1)\wndr{\wvr{R}}_{nm}^{12}(u|v|z_1 +z_2)\wndr{B}_m^{+23}(v|z_2)\non\\
&\qquad\cdot\wndr{B}_m^{23}(v|z_2+hc/2)^{-1}\wndr{\wvr{R}}_{nm}^{12}(u|v|z_1 +z_2+2hc)^{-1}\wndr{B}_n^{13}(u|z_1+hc/2)^{-1}.\label{loc4}
\end{align}

Let us  consider the second summand in \eqref{mod3}. First, by swapping tensor factors $n+m+1$ and $n+m+2$ in \eqref{loc1} we get
$$\wvr{R}_{nm}^{12}(u|v|z)^{-1}T_m^{+23}(v|0)\wvr{R}_{nm}^{12}(u|v|z-2hc)T_n^{+14}(u|0)(\vac\ot\vac).
$$
Next, by applying $Y_{\Wc_c}(z_2)(1\ot Y_{\Wc_c}(z_1))$  we obtain
\begin{align}
&\wvr{R}_{nm}^{12}(u|v|z)^{-1}\wndr{B}_m^{+23}(v|z_2)\wndr{B}_m^{23}(v|z_2+hc/2)^{-1}\non\\
&\qquad\cdot\wvr{R}_{nm}^{12}(u|v|z-2hc)
\wndr{B}_n^{+13}(u|z_1)\wndr{B}_n^{13}(u|z_1+hc/2)^{-1}.\label{loc6}
\end{align}

We  now want to apply relation  \eqref{rbrb32}   on \eqref{loc6}. However, the factors
\beq\label{loc7}
\wvr{R}_{nm}^{12}(u|v|z)^{-1}\fand \wvr{R}_{nm}^{12}(u|v|z-2hc),\qquad\text{where}\quad z=z_1-z_2,
\eeq
in \eqref{loc6} should be expanded in nonnegative powers of $z_2$, while \eqref{rbrb32}  requires for the $R$-matrices in \eqref{loc7} to be expanded in nonnegative powers of $z_1$. Fix  an integer $k\geqslant 0$. For any choice of  integers $a_1 , \ldots, a_n\geqslant 0$ and
$b_1 , \ldots , b_m\geqslant 0$ there exist an integer $r \geqslant 0$ such that the coefficients of all monomials $u_1^{a'_1}\ldots u_n^{a'_n}v_1^{b'_1}\ldots v_m^{b'_m}$,
where $0\leqslant a'_i\leqslant a_i$ and $0\leqslant b'_j\leqslant b_j$, in
\beq\label{zt3}
(z_1 -z_2)^{r }\,\wvr{R}_{nm}^{12}(u|v|z_1 -z_2)^{-1}\fand (z_1 -z_2)^{r }\,\wvr{R}_{nm}^{12}(u|v|z_1-z_2-2hc)
\eeq
coincide with the corresponding coefficients in
\beq\label{zt4}
(z_1 -z_2)^{r }\,\wvr{R}_{nm}^{12}(u|v|-z_2 +z_1)^{-1}\fand (z_1 -z_2)^{r }\,\wvr{R}_{nm}^{12}(u|v|-z_2+z_1-2hc)
\eeq
modulo $h^k$. Moreover, assume that the integer $r$ is large enough, so that the
coefficients of all monomials $u_1^{a'_1}\ldots u_n^{a'_n}v_1^{b'_1}\ldots v_m^{b'_m}$,
where $0\leqslant a'_i\leqslant a_i$ and $0\leqslant b'_j\leqslant b_j$,  in
\beq\label{zt1}
(z_1 +z_2)^{r }\,\wndr{\wvr{R}}_{nm}^{12}(u|v|z_1 +z_2)\fand (z_1 +z_2)^{r }\,\wndr{\wvr{R}}_{nm}^{12}(u|v|z_1+z_2+2hc)^{-1}
\eeq
coincide with the corresponding coefficients in
\beq\label{zt2}
(z_1 +z_2)^{r }\,\wndr{\wvr{R}}_{nm}^{12}(u|v|z_2 +z_1)\fand (z_1 +z_2)^{r }\,\wndr{\wvr{R}}_{nm}^{12}(u|v|z_2+z_1+2hc)^{-1}
\eeq
modulo $h^k$. 
By using \eqref{rbrb32} and unitarity \eqref{uni} we obtain 
\begin{align*}
&\wndr{B}_m^{23}(v|z_2+hc/2)^{-1}\wvr{R}_{nm}^{12}(u|v|-z_2+z_1-2hc)
\wndr{B}_n^{+13}(u|z_1)=\wndr{\wvr{R}}_{nm}^{12}(u|v|z_2+z_1)\\
&\qquad\cdot\wndr{B}_n^{+13}(u|z_1)\wvr{R}_{nm}^{12}(u|v|-z_2+z_1) \wndr{B}_m^{23}(v|z_2+hc/2)^{-1}
\wndr{\wvr{R}}_{nm}^{12}(u|v|z_2+z_1+2hc)^{-1}.
\end{align*}
This implies, due to the  fact that certain coefficients in \eqref{zt3} and \eqref{zt1} coincide with the corresponding coefficients in \eqref{zt4} and \eqref{zt2} modulo $h^k$, that the product of \eqref{loc6} and $(z_1^2-z_2^2)^{2r}$ coincides with
\begin{align}
&(z_1^2-z_2^2)^{2r}\, \wvr{R}_{nm}^{12}(u|v|z)^{-1}\wndr{B}_m^{+23}(v|z_2)\wndr{\wvr{R}}_{nm}^{12}(u|v|z_1+z_2)\wndr{B}_n^{+13}(u|z_1)\wvr{R}_{nm}^{12}(u|v|z)\non\\
&\qquad\qquad\cdot\wndr{B}_m^{23}(v|z_2+hc/2)^{-1}
\wndr{\wvr{R}}_{nm}^{12}(u|v|z_1+z_2+2hc)^{-1}\wndr{B}_n^{13}(u|z_1+hc/2)^{-1}\label{loc9}
\end{align}
modulo $h^k$. Finally, we rewrite \eqref{loc9} using \eqref{rbrb12}, thus getting
\begin{align}
&(z_1^2-z_2^2)^{2r}\wndr{B}_n^{+13}(u|z_1)\wndr{\wvr{R}}_{nm}^{12}(u|v|z_1+z_2)\wndr{B}_m^{+23}(v|z_2)\non\\
&\qquad\cdot\wndr{B}_m^{23}(v|z_2+hc/2)^{-1}
\wndr{\wvr{R}}_{nm}^{12}(u|v|z_1+z_2+2hc)^{-1}\wndr{B}_n^{13}(u|z_1+hc/2)^{-1}\label{loca}.
\end{align}
Since \eqref{loca}  is equal to   the product of $(z_1^2-z_2^2)^{2r}$ and \eqref{loc4}, we conclude that \eqref{mod3} holds.
\end{prf}

As with the operator $\Tc(z)=Y(T^+(0)\vac,z)$, see \cite[2.1.4]{EK}, the proof of Proposition \ref{proploc}  implies that the operator
$\Bc(z)=Y_{\Wc_c(\gl_N)}(T^+(0)\vac,z)$ satisfies the (slightly modified version of the) reflection equation  from \cite{RS}. More precisely, for any integer $n\geqslant 0$ there exist an integer $r\geqslant 0$ such that 
\begin{align*}
&(z_1^2 -z_2^2)^r\, \Bc_1(z_1) \wvr{R}_{12}(z_1-z_2 +2hc)^{-1}\Bc_2 (z_2)\wvr{R}_{12}(z_1 -z_2)\\
&\qquad\qquad \hsym\,\, (z_1^2 -z_2^2)^r\, \wvr{R}_{12}(z_1-z_2)^{-1}\Bc_2 (z_2)\wvr{R}_{12}(z_1-z_2-2hc)
\Bc_1(z_1).
\end{align*}

\section{Image of the center  \texorpdfstring{$\z(\Vc_{2c}(\gl_N))$}{z(V2c(glN))}}\label{sec3}
In this section, we  employ   map $  Y_{\Wc_c(\gl_n)}(z)$ to find explicit formulae for families of central elements in the completed algebra $\wtld{\Ac}_{c}(\gl_N)$. As a consequence, we obtain families of invariants of the vacuum module $\Wc_c(\gl_n)$. Also, we show that the image of the center  $\z(\Vc_{2c}(\gl_N))$, with respect to the map $a\mapsto Y_{\Wc_c(\gl_n)}(a,z)$, is commutative.

\subsection{Central elements of the completed algebra  \texorpdfstring{$\wtld{\Ac}_{-N/2}(\gl_N)$}{A-N/2(glN)}}

Let $I_p$ for $p\geqslant 1$ denote the left ideal of the double Yangian $\DY_c(\gl_N)$ at the level 
 $c\in\CC$, generated by all elements $t_{ij}^{(r)}$ with $r\geqslant p$.  
 As in \cite{JKMY}, define the completed  double Yangian $\wtld{\DY}_c(\gl_N)$ at the level 
  $c$    as  the $h$-adic  completion of the inverse limit
 $\lim_{\longleftarrow} \,  \DY_c(\gl_N)/  I_p.
 $
Introduce the algebra    $\wtld{\Ac}_c(\gl_N)$  as  the $h$-adic  completion of the inverse limit
$$\lim_{\longleftarrow} \,  \Ac_c(\gl_N)/ (\Ac_c(\gl_N)\cap I_p).
$$

In order to employ certain results from \cite{JKMY}, we  briefly recall the fusion procedure for the rational $R$-matrix originated in \cite{J}; see also \cite[Section 6.4]{M} for more details. Let $\mu$ be a Young diagram with $n$ boxes, whose length is less than or equal to $N$,  and let $\Uc$ be a standard $\mu$-tableau with entries $1,\ldots,n$. For  $k=1,\ldots,n$ define the contents $c_k$ of $\Uc$ by $c_k = j-i$ if $k$ occupies the box $(i,j)$ of $\Uc$.
Denote by $e_{\Uc}$   the primitive idempotent in the group algebra $\CC[\Sym_n]$ of the symmetric group $\Sym_n$, which is associated with $\Uc$
through the use of the orthonormal Young bases in the irreducible representations of $\Sym_n$.
The  group $\Sym_n$ acts on the space $(\CC^N)^{\ot\ts n}$ by permuting the tensor
factors. Denote by $\Ec^{}_{\Uc}$ the image
of $e^{}_{\Uc}$ with respect to this action.
By \cite{J}, the consecutive evaluations $u_1 =hc_1,\ldots , u_n =hc_n$ 
of the function 
$$R(u_1,\dots,u_n)\coloneqq \prod_{1\leqslant i< j\leqslant n} R_{ij}(u_i -u_j),$$
where the product is taken in the lexicographical order on the pairs $(i,j)$,
 are well-defined. Furthermore,
  the result is proportional to $\Ec^{}_{\Uc}$, i.e. 
\beq\label{fusion}
R(u_1,\dots,u_n)\big|_{u_1=hc_1}
\big|_{u_2=hc_2}\dots \big|_{u_n=hc_n}=p(\mu)\ts \Ec^{}_{\Uc},
\eeq
where $p(\mu)$ denotes the product of all hook lengths of the boxes of $\mu$.

Let
\beq\label{uovi}
u_\mu =(u_1,\ldots,u_n),\qquad\text{where}\qquad
u_k = u+hc_k \quad\text{for }k=1,\ldots,n.\eeq
It was proved in \cite{JKMY} that all  coefficients of the series
\beq\label{ln}
\TT_{\mu}^+ (u)=\tr_{1,\ldots,n} \, \Ec_{\Uc} T_1^+ (u_1)\ldots T_n^{+}(u_n) \vac\in \Vc_{-N}(\gl_N)[[u]],
\eeq
 where the trace   is taken over all $n$ copies of $\ndo\CC^N$ in \eqref{ln},
belong to the center $\z(\Vc_{-N}(\gl_N))$.
The series $\TT_{\mu}^+ (u)$ does not depend on the choice of the standard $\mu$-tableau $\Uc$; see \cite{O}.
The image of the constant term in \eqref{ln}, with respect to    map \eqref{Y},  equals
\beq\label{Yln}
Y_{\Wc_{-N/2}(\g_N)}(\TT_{\mu}^+(0),u)=
\tr_{1,\ldots,n}\, \Ec_{\Uc} \wndr{B}^+_n(u_\mu)\wndr{B}_n(u_\mu -hN/4)^{-1}
\eeq
and belongs to $\om(\Wc_{-N/2}(\gl_N),\Wc_{-N/2}(\gl_N)((u))[[h]] )$. All coefficients of  series  \eqref{Yln},
\begin{align}
&\wtld{\BB}_{\mu} (u)\coloneqq\tr_{1,\ldots,n}\, \Ec_\Uc \wndr{B}^+_n(u_\mu)\wndr{B}_n(u_\mu-hN/4)^{-1}\label{c1}
\end{align}
 can be also viewed   as elements of the completed algebra
$\wtld{\Ac}_{-N/2}(\gl_N)$.

Consider the tensor product
\beq\label{ten1ten1}
\ndo\CC^N \ot (\ndo\CC^N)^{\ot n}\ot \wtld{\Ac}_{-N/2}(\gl_N),
\eeq
where   the $n+1$ copies
of $\ndo\CC^N$ are now labeled by $0,\ldots,n$.
It will be convenient to denote the tensor factors $\ndo\CC^N$, $(\ndo\CC^N)^{\ot n}$ and $\wtld{\Ac}_{-N/2}(\gl_N)$ in \eqref{ten1ten1}  by the  superscripts $0,1$ and $2$
respectively,  so that, e.g., for the variable $u_0$ and variables \eqref{uovi}  we have
\begin{align*}
& \wvr{\wndr{R}}_{1n}^{01}(u_0|u_\mu)=  \wvr{R}_{01}(u_0+u_1)\ldots \wvr{R}_{0n}(u_0+u_n).
\end{align*}
 The arrow at the top of the symbol will indicate  that the products are written in the opposite order, e.g., 
\begin{align*}
& \cev{\wvr{\wndr{R}}}_{1n}^{01}(u_0|u_\mu)=\wvr{R}_{0n}(u_0+u_{n})\ldots \wvr{R}_{01}(u_0+u_1).
\end{align*}

\begin{lem}
The following equalities hold on $\ndo\CC^N \ot (\ndo\CC^N)^{\ot n}\ot \wtld{\DY}_{c}(\gl_N)$:
\begin{align}
&\Ec_{\Uc} \wvr{R}_{1n}^{01}(u_0|u_\mu)= \cev{\wvr{R}}_{1n}^{01}(u_0|u_\mu)\Ec_{\Uc},\qquad \Ec_{\Uc} \wvr{R}_{1n}^{01}(u_0|u_\mu)^{-1}= \cev{\wvr{R}}_{1n}^{01}(u_0|u_\mu)^{-1}\Ec_{\Uc},\label{xxx1}\\
&\Ec_{\Uc} \wvr{\wndr{R}}_{1n}^{01}(u_0|u_\mu)= \cev{\wvr{\wndr{R}}}_{1n}^{01}(u_0|u_\mu)\Ec_{\Uc},\qquad \Ec_{\Uc} \wvr{\wndr{R}}_{1n}^{01}(u_0|u_\mu)^{-1}= \cev{\wvr{\wndr{R}}}_{1n}^{01}(u_0|u_\mu)^{-1}\Ec_{\Uc},\label{xxx2}\\
&\Ec_{\Uc}\wndr{B}^{+12}_n(u_\mu)=\cev{\wndr{B}}^{+12}_n(u_\mu)\Ec_{\Uc},\qquad \Ec_{\Uc}\wndr{B}_n^{12}(u_\mu-hN/4)^{-1}=\cev{\wndr{B}}_n^{12}(u_\mu-hN/4)^{-1}\Ec_{\Uc},\label{xxx4}\\
&\Ec_{\Uc}T^{+12}_n(u_\mu|0)=\cev{T}^{+12}_n(u_\mu|0)\Ec_{\Uc},\qquad \Ec_{\Uc}\cev{T}_n^{+12}(-u_\mu|0)^{-1}=T_n^{+12}(-u_{\mu}|0)^{-1}\Ec_{\Uc},\label{xxx5}\\
&\Ec_{\Uc}\cev{T}^{12}_n(-u_\mu+hN/4|0)=T^{12}_n(-u_\mu+hN/4|0)\Ec_{\Uc} ,\label{xxx6}\\
&\Ec_{\Uc}T_n^{12}(u_\mu-3hN/4|0)^{-1}=\cev{T}_n^{12}(u_\mu-3hN/4|0)^{-1}\Ec_{\Uc},\label{xxx7}
\end{align}
where   $\Ec_{\Uc}$  is applied on the tensor factors $1,\ldots,n$, i.e. $\Ec_{\Uc}$ denotes the operator $1\ot\Ec_{\Uc}$ on $\ndo\CC^N \ot (\ndo\CC^N)^{\ot n}$.
\end{lem}

\begin{prf}
The given equalities follow from  fusion procedure \eqref{fusion} with the use of Yang--Baxter equation \eqref{ybe}, unitarity \eqref{uni}  and   relations   \eqref{RTT2}--\eqref{RTT3} and \eqref{rbrb1}--\eqref{rbrb3}. More details on the proof can be found in \cite[Proof of Theorem 2.4]{JKMY} (for relations \eqref{xxx1}--\eqref{xxx2}), first part of the proof of Theorem \ref{zmaina} (for relations \eqref{xxx4}) and in \cite[Proof of Theorem 3.2]{FJMR} (for relations \eqref{xxx5}--\eqref{xxx7}). 
As an illustration, let us prove the first equality in \eqref{xxx2}.
For the variables $v=(u+v_1,\ldots ,u+v_n)$ 
Yang--Baxter equation \eqref{ybe} implies
\beq\label{phya}\prod_{1\leqslant i< j\leqslant n} R_{ij}(v_i -v_j)\,\cdot \,
\wvr{\wndr{R}}_{1n}^{01}(u_0|v)
=
 \cev{\wvr{\wndr{R}}}_{1n}^{01}(u_0|v)\,\cdot \,
\prod_{1\leqslant i< j\leqslant n} R_{ij}(v_i -v_j),\eeq
where the products are written in the lexicographical order on the pairs $(i,j)$.
By  applying consecutive evaluations $v_1 =hc_1,\ldots ,v_n =hc_n$ on \eqref{phya} and using \eqref{fusion} we get
$
 \Ec_{\Uc} \wvr{R}_{1n}^{01}(u_0|u_\mu)= \cev{\wvr{R}}_{1n}^{01}(u_0|u_\mu)\Ec_{\Uc}, 
$
as required.
\end{prf}

The following is our main result in this section. Its proof adapts the standard $R$-matrix techniques used with $RTT$ relations, see, e.g., \cite[Theorem 3.2]{FJMR}, to the reflection algebra setting.

\begin{thm}\label{main}
All coefficients of $\wtld{\BB}_{\mu} (u)$   belong to the center of the   algebra $\wtld{\Ac}_{-N/2}(\gl_N)$.
\end{thm}

\begin{prf}
We first prove  that for the variable $u_0$ and variables \eqref{uovi} the following equality holds  on $\ndo\CC^N \ot \wtld{\Ac}_{-N/2}(\gl_N)$:
\beq\label{c2}
 B(u_0)\wtld{\BB}_{\mu} (u)=\wtld{\BB}_{\mu} (u)B(u_0).
\eeq

By applying $B_0(u_0)$ on \eqref{c1} and using   notation  as in  \eqref{ten1ten1} we get
\beq\label{c3}
\tr_{1,\ldots,n}\, \Ec_\Uc B_0(u_0)\wndr{B}^{+12}_n(u_\mu)\wndr{B}_n^{12}(u_\mu-hN/4)^{-1}.
\eeq
As with the proof of \eqref{csym5}, we employ \eqref{rbrb3} and \eqref{csym2} to express \eqref{c3} as
\begin{align}
\tr_{1,\ldots,n} \,&\Ec_\Uc\,\bigg(
{}^{rl\hspace{-3pt}}\left(\wndr{\wvr{R}}_{1n}^{01}(u_0-3hN/4|u_\mu)^{-1}\right) \cdot\Big( \wvr{R}_{1n}^{01}(u_0-3hN/4|u_\mu)^{-1}\wndr{B}_{n}^{+12}(u_\mu)\Big.\bigg.\non\\
&\bigg.\Big. \wndr{\wvr{R}}_{1n}^{01}(u_0-3hN/4|u_\mu)B_{0}(u_0)\wvr{R}_{1n}^{01}(u_0+hN/4|u_\mu) \wndr{B}_n^{12}(u_\mu-hN/4)^{-1} \Big)\bigg)
.\label{c5}
\end{align}
Since $\Ec_\Uc^2=\Ec_\Uc$, the second equality in \eqref{xxx2}  implies
\beq\label{c4}
\Ec_\Uc K= \Ec_\Uc^2 K=\Ec_\Uc \cev{K} \Ec_\Uc =\Ec_\Uc \cev{K} \Ec_\Uc^2\qquad\text{for }K=\wndr{\wvr{R}}_{1n}^{01}(u_0-3hN/4|u_\mu)^{-1}.
\eeq
By using \eqref{c4} we can write \eqref{c5} as
\begin{align*}
\tr_{1,\ldots,n}\, &\Ec_\Uc \,\bigg(
{}^{rl\hspace{-3pt}}\left(\cev{K}\right) \cdot   \Big(\Ec_\Uc ^2 \wvr{R}_{1n}^{01}(u_0-3hN/4|u_\mu)^{-1}\wndr{B}_{n}^{+12}(u_\mu)\Big.\bigg.\non\\
&\bigg.\Big. \wndr{\wvr{R}}_{1n}^{01}(u_0-3hN/4|u_\mu)B_{0}(u_0)\wvr{R}_{1n}^{01}(u_0+hN/4|u_\mu)\wndr{B}_n^{12}(u_\mu-hN/4)^{-1}\Big)
\bigg).
\end{align*}
Due to the cyclic property of the trace, this equals to
\begin{align}
\tr_{1,\ldots,n}\,& \Ec_\Uc \wvr{R}_{1n}^{01}(u_0-3hN/4|u_\mu)^{-1}\wndr{B}_{n}^{+12}(u_\mu)\wndr{\wvr{R}}_{1n}^{01}(u_0-3hN/4|u_\mu)\non\\
&B_{0}(u_0)\wvr{R}_{1n}^{01}(u_0+hN/4|u_\mu)
\wndr{B}_n^{12}(u_\mu-hN/4)^{-1}\Ec_\Uc 
\cev{K} \Ec_\Uc .\label{c8}
\end{align}
By $\Ec_\Uc^2=\Ec_\Uc$ and the second equality in \eqref{xxx1}  we have
$$\Ec_\Uc L=\Ec_\Uc^2 L=\Ec_\Uc \cev{L}\Ec_\Uc =\Ec_\Uc \cev{L}\Ec_\Uc^2 =\Ec_\Uc^2 L\Ec_\Uc =\Ec_\Uc L\Ec_\Uc \quad\text{for }
L=\wvr{R}_{1n}^{01}(u_0-3hN/4|u_\mu)^{-1}.$$
Therefore, using the cyclic property of the trace and $\Ec_\Uc^2=\Ec_\Uc$, we can write \eqref{c8} as
\begin{align*}
\tr_{1,\ldots,n}\,& \wvr{R}_{1n}^{01}(u_0-3hN/4|u_\mu)^{-1}\Ec_\Uc \wndr{B}_{n}^{+12}(u_\mu)\wndr{\wvr{R}}_{1n}^{01}(u_0-3hN/4|u_\mu)\non\\
&B_{0}(u_0)\wvr{R}_{1n}^{01}(u_0+hN/4|u_\mu)
\wndr{B}_n^{12}(u_\mu-hN/4)^{-1}\Ec_\Uc 
\cev{K} \Ec_\Uc .
\end{align*}
We now employ first equalities in \eqref{xxx1} and \eqref{xxx2}, together with \eqref{xxx4}, to move the leftmost copy of $\Ec_\Uc $ to the right, which gives us:
\begin{align}
\tr_{1,\ldots,n}\,&  \wvr{R}_{1n}^{01}(u_0-3hN/4|u_\mu)^{-1}\cev{\wndr{B}}_{n}^{+12}(u_\mu)\cev{\wndr{\wvr{R}}}_{1n}^{01}(u_0-3hN/4|u_\mu)\non\\
&B_{0}(u_0)\cev{\wvr{R}}_{1n}^{01}(u_0+hN/4|u_\mu)
\cev{\wndr{B}}_n^{12}(u_\mu-hN/4)^{-1}\Ec_\Uc^2 
\cev{K}\Ec_\Uc .\label{c6}
\end{align}
Using \eqref{c4} and $\Ec_\Uc^2=\Ec_\Uc$ we   replace $\Ec_\Uc^2 
\cev{K} \Ec_\Uc $  with $\Ec_\Uc 
K=\Ec_\Uc \wndr{\wvr{R}}_{1n}^{01}(u_0-3hN/4|u_\mu)^{-1}$ in \eqref{c6}.  
Next, we employ the first equalities in \eqref{xxx1} and \eqref{xxx2}, together with \eqref{xxx4}, to move the remaining copy of $\Ec_\Uc$ to the left, thus getting
\begin{align}
\tr_{1,\ldots,n}\,&  \wvr{R}_{1n}^{01}(u_0-3hN/4|u_\mu)^{-1}\Ec_\Uc \wndr{B}_{n}^{+12}(u_\mu)\wndr{\wvr{R}}_{1n}^{01}(u_0-3hN/4|u_\mu)\non\\
&B_{0}(u_0)\wvr{R}_{1n}^{01}(u_0+hN/4|u_\mu)
\wndr{B}_n^{12}(u_\mu-hN/4)^{-1}
\wndr{\wvr{R}}_{1n}^{01}(u_0-3hN/4|u_\mu)^{-1}.\label{c7}
\end{align}
By applying \eqref{rbrb2} on the last four factors in \eqref{c7} and then by canceling the adjacent terms $\wndr{\wvr{R}}_{1n}^{01}(u_0-3hN/4|u_\mu)^{\pm 1}$ we obtain
\begin{align*}
\tr_{1,\ldots,n}\,&  \wvr{R}_{1n}^{01}(u_0-3hN/4|u_\mu)^{-1}\Ec_\Uc \wndr{B}_{n}^{+12}(u_\mu)
\wndr{B}_n^{12}(u_\mu-hN/4)^{-1}
\wvr{R}_{1n}^{01}(u_0+hN/4|u_\mu)B_{0}(u_0).
\end{align*}

In order to prove  \eqref{c2}, it is sufficient to verify that the expression
\begin{align}\label{c9}
\tr_{1,\ldots,n}\,&  \wvr{R}_{1n}^{01}(u_0-3hN/4|u_\mu)^{-1}\Ec_\Uc \wndr{B}_{n}^{+12}(u_\mu)
\wndr{B}_n^{12}(u_\mu-hN/4)^{-1}
\wvr{R}_{1n}^{01}(u_0+hN/4|u_\mu)
\end{align}
is equal to $\wtld{\BB}_{\mu} (u)$. By the property
$\tr_{1,\ldots, n}\, XY=\tr_{1,\ldots, n}\, X^{t_1\ldots t_n}Y^{t_1\ldots t_n}$
for 
\begin{align*}
X=\wvr{R}_{1n}^{01}(u_0-3hN/4|u_\mu)^{-1}\Ec_\Uc \wndr{B}_{n}^{+12}(u_\mu)
\wndr{B}_n^{12}(u_\mu-hN/4)^{-1}\hspace{-4pt}\fand\hspace{-4pt}
Y=\wvr{R}_{1n}^{01}(u_0+hN/4|u_\mu)
\end{align*}
we conclude that \eqref{c9} is equal to
\begin{align*}
&\tr_{1,\ldots,n} \, \left(\Ec_\Uc \wndr{B}_{n}^{+12}(u_\mu)
\wndr{B}_n^{12}(u_\mu-hN/4)^{-1}\right)^{t_1\ldots t_n}  Z,\qquad
\text{where}\\
&Z=\left(\wvr{R}_{1n}^{01}(u_0-3hN/4|u_\mu)^{-1}\right)^{t_1\ldots t_n}
\wvr{R}_{1n}^{01}(u_0+hN/4|u_\mu)^{t_1\ldots t_n}.
\end{align*}
Finally,  crossing symmetry property \eqref{csym} implies
$Z=1$,
so  \eqref{c2} clearly follows.

Consider the tensor product
\beq\label{ten2}
(\ndo\CC^N)^{\ot n}\ot \ndo\CC^N \ot \wtld{\Ac}_{-N/2}(\gl_N) ,
\eeq
where the $n+1$ copies
of $\ndo\CC^N$ are now labeled by $1,\ldots,n+1$. 
It will be convenient to denote
 the tensor factors 
$ (\ndo\CC^N)^{\ot n}$, $\ndo\CC^N$ and $\wtld{\Ac}_{-N/2}(\gl_N)$ in \eqref{ten2}
 by the superscripts 1,2 and 3
respectively.\footnote{We introduce the new labeling   because the application of the original labels, as in \eqref{ten1ten1}, would require   different, more appropriate notation. For example, notice that the $R$-matrices in the first part of the proof should be expanded in nonnegative powers of the variable $u$, while the $R$-matrices in the following, second part of the proof should be expanded in nonnegative powers of the variable $u_{n+1}$.}
Our next goal is to prove that for variables \eqref{uovi} and the variable $u_{n+1}$ the following equality holds on  $\ndo\CC^N\ot \wtld{\Ac}_{-N/2}(\gl_N)$:
\beq\label{cd2}
 B^+(u_{n+1})\wtld{\BB}_{\mu} (u)=\wtld{\BB}_{\mu} (u)B^{+}(u_{n+1}).
\eeq

The proof of \eqref{cd2} is analogous to the proof of \eqref{c2}, so we only briefly sketch some details  to take care of minor differences.
First, by applying $B_{n+1}^+(u_{n+1})$ on \eqref{c1} and using notation \eqref{ten2}  we get
\beq\label{cd3}
\tr_{1,\ldots,n} \, \Ec_\Uc  B_{n+1}^+(u_{n+1})\wndr{B}^{+13}_n(u_\mu)\wndr{B}_n^{13}(u_\mu-hN/4)^{-1}.
\eeq
As with the proof of \eqref{csym5}, we employ \eqref{rbrb1} and \eqref{csym2} to express \eqref{cd3} as
\begin{align}
\tr_{1,\ldots,n}\,  &\Ec_\Uc  \bigg(
{}^{lr\hspace{-3pt}}\left(\wndr{\wvr{R}}_{n1}^{12}(u_\mu-hN|u_{n+1})^{-1}\right)\cdot\Big( \wvr{R}_{n1}^{12}(u_\mu|u_{n+1})\wndr{B}_{n}^{+13}(u_\mu)\Big.\bigg.\non\\
&\Big.\bigg. \wndr{\wvr{R}}_{n1}^{12}(u_\mu|u_{n+1}) B_{n+1}^+(u_{n+1})\wvr{R}_{n1}^{12}(u_\mu|u_{n+1})^{-1}
\wndr{B}_n^{13}(u_\mu-hN/4)^{-1}\Big)\bigg)
 .\label{cd5}
\end{align}
We may now proceed as in the first part of the proof and, using the properties of the primitive idempotent $\Ec_\Uc$, show that \eqref{cd5} is equal to
\begin{align}
\tr_{1,\ldots,n} \,&
 \wvr{R}_{n1}^{12}(u_\mu|u_{n+1})\Ec_\Uc \wndr{B}_{n}^{+13}(u_\mu)\wndr{\wvr{R}}_{n1}^{12}(u_\mu|u_{n+1})\non\\
&B_{n+1}^+(u_{n+1})\wvr{R}_{n1}^{12}(u_\mu|u_{n+1})^{-1}
\wndr{B}_n^{13}(u_\mu-hN/4)^{-1}\wndr{\wvr{R}}_{n1}^{12}(u_\mu-hN|u_{n+1})^{-1}.\label{cd6}
\end{align}
By applying \eqref{rbrb3} to the last four factors in \eqref{cd6} and then canceling the adjacent terms $\wndr{\wvr{R}}_{n1}^{12}(u_\mu|u_{n+1})^{\pm 1}$ we obtain
\begin{align*}
\tr_{1,\ldots,n}\, &
 \wvr{R}_{n1}^{12}(u_\mu|u_{n+1})\Ec_\Uc \wndr{B}_{n}^{+13}(u_\mu)
 \wndr{B}_n^{13}(u_\mu-hN/4)^{-1}
 \wvr{R}_{n1}^{12}(u_\mu-hN|u_{n+1})^{-1}
B_{n+1}^+(u_{n+1}).
\end{align*}
Finally, in order to prove \eqref{cd2}, it is sufficient to check that the expression
\begin{align*} 
\tr_{1,\ldots,n} \,&
 \wvr{R}_{n1}^{12}(u_\mu|u_{n+1})\Ec_\Uc \wndr{B}_{n}^{+13}(u_\mu)
 \wndr{B}_n^{13}(u_\mu-hN/4)^{-1}
 \wvr{R}_{n1}^{12}(u_\mu-hN|u_{n+1})^{-1}
\end{align*}
is equal to $\wtld{\BB}_{\mu} (u)$. This can be done as in the first part of the proof, by employing 
crossing symmetry property \eqref{csym} and unitarity \eqref{uni}.

The statement of the theorem now follows from  \eqref{c2} and \eqref{cd2}.
\end{prf}

We now consider two   special cases of Theorem \ref{main}. 
Denote by $\mu^{\text{row}}_n$ and $\mu^{\text{col}}_n$ the row diagram with $n$ boxes and the column diagram with $n$ boxes respectively.
The unique idempotents corresponding to the  standard $\mu^{\text{row}}_n$-tableau  and $\mu^{\text{col}}_n$-tableau coincide with the images $H^{(n)}$ and $A^{(n)}$
of the symmetrizer and the anti-symmetrizer
$$
h^{(n)}=\frac{1}{n!}\sum_{s\in\Sym_n} s\Fand
a^{(n)}=\frac{1}{n!}\sum_{s\in\Sym_n} \sgn s\cdot s
$$
under the action of the symmetric group $\Sym_n$ on $(\CC^N)^{\ot n}$.
In this two cases, \eqref{c1} becomes
\begin{align*}
&\wtld{\BB}_{\mu^{\text{row}}_n} (u)=\tr_{1,\ldots,n}\, H^{(n)} \wndr{B}^+_n(u_+)\wndr{B}_n(u_+-hN/4)^{-1},\\
&\wtld{\BB}_{\mu^{\text{col}}_n} (u)=\tr_{1,\ldots,n} \, A^{(n)} \wndr{B}^+_n(u_-)\wndr{B}_n(u_- -hN/4)^{-1},
\end{align*}
where 
\beq\label{uminus}
u_\pm =(u,u\pm h,\ldots ,u\pm (n-1)h).
\eeq 
 Note that $u_+=u_{\mu^{\text{row}}_n}$ and $u_-=u_{\mu^{\text{col}}_n}$; recall \eqref{uovi}. 
Consider the  series 
\begin{align*}
&\wtld{\XX}_{\mu^{\text{row}}_n } (u)=\tr_{1,\ldots,n}\, H^{(n)}\, T^+_n(u_+|0)\,\cev{T}^+_n(-u_+|0)^{-1}\,\cev{T}_n(-u_+ +hN/4|0)\,T_n(u_+ -3hN/4|0)^{-1}, \\
&\wtld{\XX}_{\mu^{\text{col}}_n } (u)=\tr_{1,\ldots,n}\, A^{(n)}\, T^+_n(u_-|0)\,\cev{T}^+_n(-u_-|0)^{-1}\,\cev{T}_n(-u_- +hN/4|0)\,T_n(u_- -3hN/4|0)^{-1}
\end{align*}
in $\wtld{\DY}_{-N/2}(\gl_N)[[u^{\pm 1}]]$,
where,
as before, the arrows  indicate that the products are written in the opposite order, e.g., for $w=(w_1,\ldots, w_n)$ we have 
$\cev{T}_n(w|0)=T_n(w_n)\ldots T_1(w_1)$.

\begin{kor}\label{XXXxxx}
Suppose that the matrix $G$, given by \eqref{diagonal}, is equal
to $\pm I$. Then all coefficients of $\wtld{\XX}_{\mu^{\text{row}}_n }  (u)$  and $\wtld{\XX}_{\mu^{\text{col}}_n }  (u)$  belong to the center of the  algebra $\wtld{\Ac}_{-N/2}(\gl_N)$.
\end{kor}

\begin{prf}
Let $G=\varepsilon I$ for $\varepsilon\in\left\{\pm 1\right\}$.
For the family of variables $w=(w_1,\ldots, w_n)$ we have
\begin{align*}
&\wndr{B}_{n}^{+}(w)= \varepsilon^{n}\, T_n^+(w|0)  \left(
\prod_{1\leqslant i<j\leqslant n} \wvr{R}_{ij}(w_i+w_j)\right)
\cev{T}_n^+(-w|0)^{-1} \Fand\\
&\wndr{B}_{n}(w-hN/4)= \varepsilon^{n}\, T_n(w-3hN/4|0)  \hspace{-1pt}\left(
\prod_{1\leqslant i<j\leqslant n} \wvr{R}_{ij}(w_i+w_j-hN)\hspace{-3pt}\right)\hspace{-1pt}
\cev{T}_n(-w+hN/4|0)^{-1},
\end{align*}
where the products are taken in the lexicographical order on the pairs $(i,j)$.
Indeed, this easily follows from \eqref{RTT2} and \eqref{RTT1}.
Next, note that for any $1\leqslant i<j\leqslant n$ there exist functions $f_{H^{(n)}}(z)$ and $f_{A^{(n)}}(z)$ in $\CC[z^{-1}][[h]]$ satisfying
$$
H^{(n)} \wvr{R}_{ij}(z) =f_{H^{(n)}}(z) H^{(n)} \Fand A^{(n)} \wvr{R}_{ij}(z) =f_{A^{(n)}}(z) A^{(n)}.
$$
Indeed, this follows from the form of Yang $R$-matrix \eqref{yang} and the fact that for any transposition $p\in\Sym_n$ we have $p h^{(n)}=h^{(n)}$ and $p a^{(n)}=\pm  a^{(n)}$.

By combining these observations with fusion procedure \eqref{fusion} and  equalities in \eqref{xxx5}--\eqref{xxx7},  we conclude that there exist functions $\theta_{n}^{\text{row}}(z)$ and $\theta_{n}^{\text{col}}(z)$ in $\CC[z^{-1}][[h]]$ such that
\beq\label{56zt7}
\wtld{\BB}_{\mu^{\text{row}}_n } (u)=\theta_{n}^{\text{row}}(u)\wtld{\XX}_{\mu^{\text{row}}_n }  (u)\Fand \wtld{\BB}_{\mu^{\text{col}}_n } (u)=\theta_{n}^{\text{col}}(u)\wtld{\XX}_{\mu^{\text{col}}_n } (u).
\eeq
Therefore, all coefficients of $\wtld{\XX}_{\mu^{\text{row}}_n } (u)$  and $\wtld{\XX}_{\mu^{\text{col}}_n } (u)$  belong to  the  algebra $\wtld{\Ac}_{-N/2}(\gl_N)$. Finally, \eqref{56zt7} and Theorem \ref{main} imply the statement of the corollary. 
\end{prf}

It is worth noting that the functions $\theta_{n}^{\text{row}}(z)$ and $\theta_{n}^{\text{col}}(z)$ can be computed explicitly, in terms of the function  $g(u)\in 1+u^{-1}\CC[[u^{-1}]]$ defined by \eqref{geq}; cf. \cite[Theorem 3.4]{MR}.

\subsection{Invariants of the vacuum module \texorpdfstring{$\Wc_{-N/2}(\gl_N)$}{W-N/2(glN)}}
In  this section we present some further consequences of Theorem \ref{main}. Let $c$ be an arbitrary complex number.
We can view $\Wc_{c}(\gl_N)$ as a module for the  algebra $\wtld{\Ac}_{c}(\gl_N)$. Recall \eqref{diagonal} and define the {\em submodule of invariants} of $\Wc_{c}(\gl_N)$
by
$$\z(\Wc_{c}(\gl_N))=\left\{w\in\Wc_{c}(\gl_N)\,:\,  B(u)w=Gw  \right\}.$$
Clearly, 
an element $w\in \Wc_{c}(\gl_N)$ belongs to $\z(\Wc_{c}(\gl_N))$ if and only if
$$b_{ij}(u) w =\delta_{ij}\varepsilon_i w\qquad\text{for all }i,j=1,\ldots ,N,\, r=1,2,\ldots .$$
In particular,   \eqref{invs} implies that $\vac$ is an element  of  $\z(\Wc_{c}(\gl_N))$. 
Consider the series  
\begin{align}
&\BB_{\mu} (u)\coloneqq \wtld{\BB}_{\mu}(u)\vac\,\in\, \Wc_{-N/2}(\gl_N) [[u^{\pm 1}]].\label{bxn}
\end{align}
Denote by $\wht{\B}^+ (\gl_N)$ the $h$-adic completion of the algebra $\B^+(\gl_N)$. All coefficients of $\BB_{\mu} (u)$ can be viewed as elements of  $\wht{\B}^+ (\gl_N)$.
\begin{kor}\label{Bns}
All coefficients of the series $\BB_{\mu} (u)$  belong to the  submodule of invariants $\z(\Wc_{-N/2}(\gl_N))$. All coefficients  of 
$\BB_{\mu} (u)\in \wht{\B}^+ (\gl_N)[[u^{\pm 1}]]$  pairwise commute.
\end{kor}
\begin{prf}
Using Theorem \ref{main} and \eqref{bxn} we get 
\begin{align*}
b_{ij}(v)\BB_{\mu} (u)&=b_{ij}(v)\wtld{\BB}_{\mu} (u) \vac=\wtld{\BB}_{\mu} (u)b_{ij}(v) \vac=\wtld{\BB}_{\mu} (u)\delta_{ij}\varepsilon_i \vac=\delta_{ij}\varepsilon_i\wtld{\BB}_{\mu} (u) \vac=\delta_{ij}\varepsilon_i\BB_{\mu} (u)
\end{align*}
for any $i,j=1,\ldots ,N$,
which proves the first part of the corollary.

Let $\mu$ and $\nu$ be   any two partitions   having at most $N$ parts. Using Theorem \ref{main} we get
\beq\label{xyz}\wtld{\BB}_{\mu} (u)\wtld{\BB}_{\nu} (v)\vac=\wtld{\BB}_{\mu} (u)\BB_{\nu} (v)\vac=\BB_{\nu} (v)\wtld{\BB}_{\mu} (u)\vac=\BB_{\nu} (v)\BB_{\mu} (u).\eeq
Since all coefficients of the series $\wtld{\BB}_{\mu} (u)$ and $\wtld{\BB}_{\nu} (v)$  commute, we can prove analogously that
$\wtld{\BB}_{\mu} (u)\wtld{\BB}_{\nu} (v)\vac=\BB_{\mu} (u)\BB_{\nu} (v),$ 
which, together with \eqref{xyz}, implies
$[\BB_{\mu} (u),\BB_{\nu} (v)]=0,$
as required.
\end{prf}

  Corollaries \ref{XXXxxx} and \ref{Bns} imply
\begin{kor}\label{Bns2}
Let $G=\pm I$. All coefficients of the Taylor series 
\begin{align*}
&\tr_{1,\ldots,n}\, H^{(n)}\, T^+_n(u_+|0)\,\cev{T}^+_n(-u_+|0)^{-1}\vac\Fand \tr_{1,\ldots,n}\, A^{(n)}\, T^+_n(u_-|0)\,\cev{T}^+_n(-u_-|0)^{-1}\vac
\end{align*}
belong to the submodule of invariants $\z(\Wc_{-N/2}(\gl_N))$. 
\end{kor}

For any two partitions $\mu$ and $\nu$  which have at most $N$ parts we have
\beq\label{kom1}
[\wtld{\BB}_{\mu} (u),\wtld{\BB}_{\nu} (v)]=0
\eeq
in the algebra
$\wtld{\Ac}_{-N/2}(\gl_N)$.
Clearly, \eqref{kom1} remains true if we view $\wtld{\BB}_{\mu} (u)$ and $\wtld{\BB}_{\nu}(v)$ as operators on 
$\Wc_{-N/2}(\gl_N)$. Applying the substitutions $u\leftrightarrow z_1+u$ and $v\leftrightarrow z_2+v$ 
we get
\beq\label{kom2}
[\wtld{\BB}_{\mu} (z_1+u),\wtld{\BB}_{\nu} (z_2+v)]=0\quad\text{on }\Wc_{-N/2}(\gl_N).
\eeq
Note that \eqref{kom2} can be written as
\beq\label{kom3}
[Y_{\Wc_{-N/2}(\gl_N)}(\TT^{+}_{\mu} (u) ,z_1),Y_{\Wc_{-N/2}(\gl_N)}(\TT^{+}_{\nu} (v),z_2)]=0.
\eeq

\begin{thm}\label{main2}
Let $a$ be an element of the center $\z(\Vc_{-N}(\gl_N))$.
\begin{enumerate}[(1)]
\item\label{nnn111}
For any $b\in\z(\Vc_{-N}(\gl_N))$   we have
\beq\label{kom4}
[Y_{\Wc_{-N/2}(\gl_N)}(a,z_1),Y_{\Wc_{-N/2}(\gl_N)}(b,z_2)]=0 .
\eeq
\item\label{mmm222}
For any $x\in \wtld{\Ac}_{-N/2}(\gl_N)$ 
\beq\label{rrrkom4}
[Y_{\Wc_{-N/2}(\gl_N)}(a,z),x]=0\quad\text{on }\Wc_{-N/2}(\gl_N).
\eeq
\end{enumerate}
\end{thm}

\begin{prf} \eqref{nnn111}
Due to \cite[Theorem 4.9]{JKMY},  the center $\z(\Vc_{-N}(\gl_N))$ is a commutative associative algebra with respect to the product given by $a\cdot b=a_{-1}b$ for $a,b\in \z(\Vc_{-N}(\gl_N))$. Furthermore, it was proved therein that the algebra $\z(\Vc_{-N}(\gl_N))$
is topologically generated (with respect to the $h$-adic topology)  by the elements
$\Phi_m^{(r)}$, where $m=1,\ldots,N$ and $r=0,1,\ldots $, defined by
$$
\sum_{r=0}^{\infty} \Phi_m^{(r)}u^r \coloneqq h^{-m}\sum_{k=0}^m (-1)^k \binom{N-k}{m-k}\TT^{+}_{\mu_k^{\text{col}}} (u) \,\in\, \z(\Vc_{-N}(\gl_N))[[u]].
$$
By \eqref{kom3} we conclude that  \eqref{kom4} holds for any two elements $a$ and $b$ which belong to the family $\Phi_m^{(r)}$, $m=1,\ldots,N$, $r=0,1,\ldots $. 
Finally, part \eqref{aaaqqq} of Proposition \ref{protech}  implies that \eqref{kom4} holds for any two elements $a,b\in\z(\Vc_{-N}(\gl_N))$.

\eqref{mmm222} It suffices to observe that, due to Theorem \ref{main}, Equality \eqref{rrrkom4} holds if  $a=\Phi_m^{(r)}$ for some $m=1,\ldots,N$ and $r=0,1,\ldots $. 
Hence, part \eqref{bbbqqq} of  Proposition \ref{protech}  implies that \eqref{rrrkom4} holds for any   $a\in\z(\Vc_{-N}(\gl_N))$.
\end{prf}

Due to Theorem \ref{main2},  for any $a\in\z(\Vc_{-N}(\gl_N))$ and $i,j=1,\ldots , N$ we have
\beq\label{construct}
[b_{ij}(u),Y_{\Wc_{-N/2}(\gl_N)}(a,z)]=0\quad\text{on }\Wc_{-N/2}(\gl_N).
\eeq
Hence,  we can construct  elements of  $\z(\Wc_{-N/2}(\gl_N))$ as follows:
\begin{kor}
For any $a\in\z(\Vc_{-N}(\gl_N))$ and $w\in \z(\Wc_{-N/2}(\gl_N))$ all coefficients of the series $Y_{\Wc_{-N/2}(\gl_N)}(a,z)w$ belong to the submodule of invariants  $\z(\Wc_{-N/2}(\gl_N))$.
\end{kor}
\begin{prf}
Set $\Wc_{-N/2}=\Wc_{-N/2}(\g_N)$.
By employing \eqref{construct}   we get
\begin{align*}
b_{ij}(u)Y_{\Wc_{-N/2}}(a,z)w&=Y_{\Wc_{-N/2}}(a,z)b_{ij}(u)w=Y_{\Wc_{-N/2}}(a,z)\delta_{ij}\varepsilon_i w=\delta_{ij}\varepsilon_i Y_{\Wc_{-N/2}}(a,z)w
\end{align*}
 for any  $i,j=1,\ldots ,N$ and $w\in \z(\Wc_{-N/2}(\gl_N))$, 
  as required. 
\end{prf}

\subsection{Central elements and invariants at the noncritical level}
Let $c\neq -N/2$ be an arbitrary complex number.
It is well known that all coefficients of {\em quantum determinants}
\begin{align}
&\qdet T^+(u)=\sum_{\sigma\in\Sym_N}\sgn \sigma\cdot
t^+_{ \sigma(1)\tss 1}(u)\dots t^+_{ \sigma(N)\tss N}(u-(N-1)h)\,\in\, \wht{\Y}^+ (\gl_N)[[u]] ,
\label{qdetplus}\\
&\qdet T(u)=\sum_{\sigma\in\Sym_N}\sgn \sigma\cdot
t_{ \sigma(1)\tss 1}(u)\dots t_{ \sigma(N)\tss N}(u-(N-1)h)\,\in\, \Y(\gl_N) [[u^{-1}]]\label{qdet}
\end{align}
belong to the center of the algebra $\wtld{\DY}_{2c}(\gl_N)$; see, e.g., \cite[Proposition 2.8]{JKMY}.
Furthermore, if we identify the $\CC[[h]]$-modules $\wht{\Y}^+ (\gl_N)$ and $\Vc_{2c}(\gl_N)$, then all coefficients of $\qdet T^+(u)$ belong to the center $\z(\Vc_{2c}(\gl_N))$ of the quantum VOA $\Vc_{2c}(\gl_N)$; see \cite[Proposition 4.10]{JKMY}.

Set $n=N$ in \eqref{uminus}. The following equations in $(\ndo\CC^N)^{\ot N} \ot \wtld{\DY}_{2c}(\gl_N) [[u^{\pm 1}]]$ hold:
\begin{align}
A^{(N)} T_N^{+}(u_-|0)=A^{(N)}\qdet T^+(u)\fand
A^{(N)}T_N(u_-|0)=A^{(N)}\qdet T(u),\label{qdetuseful}
\end{align}
see  \cite[Section 1]{M} for more details.
By applying   quasi module map \eqref{Y} on the constant term of \eqref{qdetplus}, which is viewed as an element of the quantum VOA $\Vc_{2c}(\gl_N)$,  and by employing the first equality in \eqref{qdetuseful}, we obtain
\begin{align}\label{Ynon}
Y_{\Wc_{c}(\gl_N)}(\qdet T^+(0),u)=
\tr_{1,\ldots,N} A^{(N)}\wndr{B}^+_N(u_-)\wndr{B}_N(u_-+hc/2)^{-1}.
\end{align}
Clearly, \eqref{Ynon} belongs to $\om(\Wc_c(\gl_N),\Wc_c(\gl_N)((u))[[h]])$. However, we can view all coefficients of the series
\beq\label{cccc1}
\wtld{\BB}_{c} (u)\coloneqq \tr_{1,\ldots,N} A^{(N)}\wndr{B}^+_N(u_-)\wndr{B}_N(u_- +hc/2)^{-1}
\eeq
as elements of the  algebra
$\wtld{\Ac}_{c}(\gl_N)$, so that $\wtld{\BB}_{c} (u)$ is an element of $\wtld{\Ac}_{c}(\gl_N)[[u^{\pm 1}]]$.

\begin{pro}\label{mainc}
Let $c\neq -N/2$ be an arbitrary complex number.
\begin{enumerate}[(i)]
\item\label{cv1}
All coefficients of $\wtld{\BB}_{c} (u)$     belong to the center of the algebra $\wtld{\Ac}_{c}(\gl_N)$.
\item\label{cv2} For any $a,b\in\z(\Vc_{2c}(\gl_N))$  we have
$$[Y_{\Wc_{c}(\gl_N)}(a,z_1),Y_{\Wc_{c}(\gl_N)}(b,z_2)]=0.$$
\item\label{cv3} For any $a\in\z(\Vc_{2c}(\gl_N))$ and $x\in \wtld{\Ac}_{c}(\gl_N)$ we have
$$
[Y_{\Wc_{c}(\gl_N)}(a,z),x]=0\quad\text{on }\Wc_{c}(\gl_N).
$$
\item\label{cv4} For any $a\in\z(\Vc_{2c}(\gl_N))$ and $w\in \z(\Wc_c(\gl_N))$ all coefficients of  $Y_{\Wc_c(\gl_N)}(a,z)w$ belong to the submodule of invariants 
$\z(\Wc_c(\gl_N))$. 
\end{enumerate}
\end{pro}

\begin{prf}
\eqref{cv1} Recall that $u_-=(u_1,\ldots, u_N)=(u,\ldots, u-(N-1)h)$, so
the first equality in \eqref{qdetuseful} can be written as
\beq\label{qdetagain}
A^{(N)} T_1^{+}(u_1)\ldots T_N^{+}(u_N) =A^{(N)}\qdet T^+(u).
\eeq
We now proceed as follows (cf. \cite[Theorem 3.4]{MR}):
\begin{itemize} 
\item[$\circ$] Multiply \eqref{qdetagain} from the right by
$T_N^{+}(u_N)^{-1}\ldots T_1^{+}(u_1)^{-1} (\qdet T^+(u))^{-1}$;
\item[$\circ$] Replace $u$ with $-u+(N-1)h$;
\item[$\circ$] Conjugate the resulting equality by the permutation $(1,\ldots, N)\mapsto (N,\ldots ,1)$.
\end{itemize}
This gives us
\beq\label{qdetagain2}
A^{(N)}\qdet T^+(-u+(N-1)h)^{-1}=A^{(N)} T_1^{+}(-u_1)^{-1}\ldots T_N^{+}(-u_N)^{-1}.
\eeq
Starting from the second equality in \eqref{qdetuseful}, one can similarly prove  
\beq\label{qdetagain3}
A^{(N)}\qdet T(-u+(N-1)h-hc/2)=A^{(N)} T_N(-u_N-hc/2)\ldots T_1(-u_1-hc/2).
\eeq

By employing \eqref{qdetuseful}, \eqref{qdetagain2} and \eqref{qdetagain3} and  arguing as in the proof of Corollary \ref{XXXxxx}, we can express $\wtld{\BB}_{c} (u)$ as
\begin{align}
\wtld{\BB}_{c} (u)=\,\,&\theta_{c}(u)\,\qdet T^+ (u) \, (\qdet T^+(-u+(N-1)h))^{-1}\non\\
&\cdot\,\qdet T(-u+(N-1)h-hc/2)\, (\qdet T(u+3hc/2))^{-1}\label{qdetag}
\end{align}
for some function $\theta_c (z)$ in $\CC[z^{-1}][[h]]$.\footnote{Observe that, in contrast with the proof of 
 Corollary \ref{XXXxxx}, we no longer need the assumption $G=\pm I$  because the image of the anti-symmetrizer $A^{(N)}$ on $(\CC^N)^{\ot N}$ is one-dimensional.} Since the coefficients of quantum determinants belong to the center of the double Yangian, we conclude by   \eqref{qdetag} that the coefficients of $\wtld{\BB}_{c} (u)$ belong to the center of the algebra $\wtld{\Ac}_{c}(\gl_N)$.

\eqref{cv2}--\eqref{cv4}
By  \cite[Proposition 4.10]{JKMY}, the center $\z(\Vc_{2c}(\gl_N))$ is a commutative associative algebra with respect to the product given by $a\cdot b=a_{-1}b$ for $a,b\in \z(\Vc_{2c}(\gl_N))$. Furthermore,  it was proved therein, that the  algebra $\z(\Vc_{2c}(\gl_N))$   is topologically generated (with respect to the $h$-adic topology) by the   elements $d_0, d_1,\ldots$, which are defined    by
$$\qdet T^+ (u) =1-h(d_0 +d_1 u +d_2 u^2+\ldots).$$
Therefore,  statements \eqref{cv2}--\eqref{cv4} can be verified   using Proposition \ref{protech}, in the same way as their critical level counterparts.
\end{prf}

Consider the series 
\begin{align*}
&S^+(u)=\qdet T^+ (u)\, (\qdet T^+(-u+(N-1)h))^{-1}\,\in\, \wht{\B}^+ (\gl_N)[[u]],\\
&S^{(c)}(u)=\qdet T (u+hc) \,(\qdet T(-u+(N-1)h))^{-1}\,\in\, \B_c (\gl_N)[[u^{-1}]].
\end{align*}
By part \eqref{cv1} of Proposition \ref{mainc} and \eqref{qdetag}, all coefficients of
$S^+(u)S^{(c)}(u+hc/2)^{-1}$
belong to the center of the algebra $\wtld{\Ac}_{c}(\gl_N)$. Moreover, by applying the given expression on $\vac\in \Wc_c(\gl_N)$ and employing part \eqref{cv4} of Proposition \ref{mainc}, we see that all coefficients  of the series  $S^+(u)\vac$ belong to the submodule of invariants $\z(\Wc_c(\gl_N))$.

\begin{rem}
Let $h=1$ and $c=0$. The series $S^{(0)} (u)$    coincides, modulo the multiplicative factor from $\CC(u)$, with the {\em Sklyanin determinant} $\sdet B(u)$, whose odd coefficients are algebraically independent  and generate the center of the reflection algebra    
$\mathcal{B}(N,N-M)$; see  \cite[Theorem 3.4]{MR}.
\end{rem}

\section*{Acknowledgement}
The author would like to thank Alexander Molev for useful discussions.
The research was partially supported by the Australian Research Council and by the Croatian Science Foundation under the project 2634.

\end{document}